\documentclass[review,sort&compress,1p,times]{elsarticle}

\usepackage{amsmath}
\usepackage{pdfpages}
\usepackage{geometry}  
\usepackage[colorlinks,linkcolor=green]{hyperref}   
\usepackage{float}
\usepackage{amsfonts,amssymb}
\usepackage{bm}  
\usepackage{graphicx}
\usepackage{mathrsfs}
\usepackage{cancel} 
\usepackage{extarrows} 
\usepackage{subeqnarray}
\usepackage{cases}   
\usepackage{xcolor} 
\usepackage{colortbl,booktabs} 
\usepackage{multirow}
\usepackage{threeparttable}  
\usepackage{subfigure}
\usepackage{caption}  
\usepackage{listings}
\usepackage{setspace}
\usepackage{lineno,hyperref}
\newcommand{\ud}{\,\mathrm{d}}

\begin{document}
\begin{frontmatter}
\title{Stabilization-free virtual element method for 2D third medium contact}  

\author[mymainaddress]{Bing-Bing Xu\corref{mycorrespondingauthor}}
\cortext[mycorrespondingauthor]{Corresponding author}
\ead{bingbing.xu@ikm.uni-hannover.de}

\author[mymainaddress]{Peter Wriggers}
\address[mymainaddress]{Institute of Continuum Mechanics, Leibniz University Hannover, Hannover, Germany}

\begin{abstract}
The third medium contact has been proven to be an effective approach for simulating contact problems involving large deformations. 
Unlike traditional contact algorithms, 
the third medium contact introduces a third medium between two contacting bodies, 
thereby avoiding the complex treatment of the contact constraints.
The approach has been successfully applied in different applications in the framework of the finite element method (FEM).
As a generalization of the finite element method, 
the virtual element method (VEM) can handle arbitrary polygonal elements,
providing greater flexibility for modeling third medium contact.
However, due to the introduction of the projection operator, 
VEM requires additional stabilization terms to control the rank of the stiffness matrix.
Moreover, the regularization term in the third medium contact formulation requires a second-order numerical scheme,
which further complicates the application of VEM to such problems.
In this work, the stabilization-free virtual element method (SFVEM) is developed and applied to solve the third medium contact problems.
Different from the traditional VEM, 
SFVEM does not require additional stabilization terms, 
which simplifies the construction of the regularization term in third medium contact.
Building upon the traditional second-order FEM framework,
we present the specific format of SFVEM for solving third medium contact, 
including the construction of high-order projection operator and the tangent stiffness matrix.
Numerical examples are provided to demonstrate the effectiveness and applicability of SFVEM in solving complex 2D third medium contact problems.
\end{abstract}

\end{frontmatter}

\section{Introduction}
\label{p10.s1}
The numerical solution technology of contact problems has a wide range of applications in engineering.
In classical computational contact mechanics,
approaches to discretizing contact include node-to-node contact, node-to-segment contact, and mortar methods \cite{Wriggers2006,Mortar1}.
The most common strategies for enforcing contact constraints are the Lagrange multiplier method and the penalty method \cite{Wriggers2006,EDMcontact}.
While the Lagrange multiplier method ensures exact constraint satisfaction, it introduces additional unknowns.
The penalty method avoids extra degrees of freedom but depends strongly on the penalty parameter, which can cause either poor accuracy or ill-conditioned stiffness matrices.
These drawbacks limit the robustness of traditional contact algorithms, especially in problems with large deformations.

To address these challenges, the third medium contact method was first proposed in \cite{TMCWriggers1} based on anisotropic modeling near the contact area.
By introducing a third medium between the contacting bodies, the contact problem is reformulated as a standard finite deformation problem.
The format was originally implemented within the FEM framework and later extended to isogeometric analysis \cite{TMCcon3} and meshless methods \cite{TMCcon1}.
A similar approach named the contact domain method was proposed in \cite{CDM1, CDM2} for large deformation frictional contact problems,
though it still requires a global contact search.

A recent breakthrough in topology optimization has advanced the third medium contact method.
As introduced in \cite{TMC3, TMC4, TMC5}, the third medium is modeled as a highly compliant material,
and the regularization is achieved by introducing gradients of the deformation gradient into the strain energy density function.
Various forms of regularization have been proposed and discussed in \cite{TMC1, TMC5}.
Furthermore, frictional contact was addressed in \cite{TMC6} through an analogy with crystal plasticity.
Since the regularization term involves higher-order derivatives, 
the finite elements with higher-order shape functions are required and the efficiency is reduced.
Inspired by the method used for gradient enhancement of continuum damage models \cite{FirstOrder},
a new regularization technique was introduced by Wriggers et al. \cite{TMC2}
based on the first-order finite elements and successfully extended to three-dimensional problems \cite{TMC2,TMCWriggers2}.
Based on the framework of FEM, 
different types of elements such as the triangular elements and tetrahedral elements have been employed to solve third medium contact.

Similar with FEM,
the virtual element method (VEM) is a generalization formualtion that can handle arbitrary polygonal or polyhedral elements \cite{Veiga2012, Veiga2014}.
The flexibility in meshing makes it easier for VEM to handle problems such as crack propagation, adaptive meshes, and contact mechanics.
In addition, compared with other polygon-based methods \cite{PFEM1, PFEM2}, 
VEM is easier to construct high-order formualtion, see \cite{VEM3D1,VEM3D2,HighXu1}.
Up to now, VEM has been successfully applied in various fields,
including linear elastic problems \cite{VEMelastic1, VEMelastic2, VEMelastic3, VEMelastic4, VEMengineer}
and nonlinear problems \cite{VEMhyper1,VEMhyper2,VEMhyper3,VEMhyper4,VEMplastic3,VEMplastic2}.
In the context of contact mechanics,
VEM has been employed to solve different contact problems \cite{VEMcontact1, VEMcontact2, VEMcontact3, VEMcontact4}.
Approaches such as the node-to-node contact or node-to-segment contact are used to discretize the contact interface.
Some review and comparison of different contact algorithms based on VEM can be found in \cite{VEMbook1}.
Up to now, VEM has not been applied to solve third medium contact problems.

In this work, we develop a virtual element formualtion to solve the 2D third medium contact problems.
The main challenge in applying VEM to third medium contact lies in the construction of the regularization term.
Due to the introduction of the projection operator in VEM,
the stabilization term is usually required to ensure the rank of the stiffness matrix.
However, the regularization term in third medium contact already involves higher-order derivatives,
which complicates the construction of the stabilization term.
To overcome this difficulty, we employ the stabilization-free virtual element method (SFVEM) for the third medium contact.
The SFVEM was first proposed in \cite{StabVEM1} and has been successfully used in different fields \cite{StabVEM2, StabVEM3, StabVEM4,StabXu1,StabXu2,StabXu3}. 
The existence of the gradient projection operator makes it easier for SFVEM to calculate higher-order derivatives of displacement, 
such as the gradient of strain gradient.
We will discuss the specific format of SFVEM for solving third medium contact problems,
including the construction of high-order projection operators and the tangent stiffness matrix.
Besides, different regularization terms for the third medium contact will be compared in the numerical examples.

Based on the above description, we can divide the paper into the following parts.
In section \ref{p10.s2}, the basic idea of the third medium contact is given.
Besides in this section, the regularization terms for the third medium contact are discussed in detail.
Then, the linearization and tangent stiffness matrix for the third medium contact are presented in section \ref{p10.s3}.
The stabilization-free virtual element method is introduced in section \ref{p10.s4},
and the specific format of SFVEM for solving third medium contact is also discussed.
In section \ref{p10.s5}, several numerical examples are provided to demonstrate the effectiveness and applicability of SFVEM in solving complex 2D third medium contact problems.
Finally, some conclusions are drawn in section \ref{p10.s6}.
Some automatic differentiation codes for the linearization are also given in the appendix.

\section{Governing equations for contact mechanics}
\label{p10.s2}
\subsection{Continuum model}
\label{p10.s2.1}
\begin{figure}[htbp]
\centering
\includegraphics[width=1.0\textwidth]{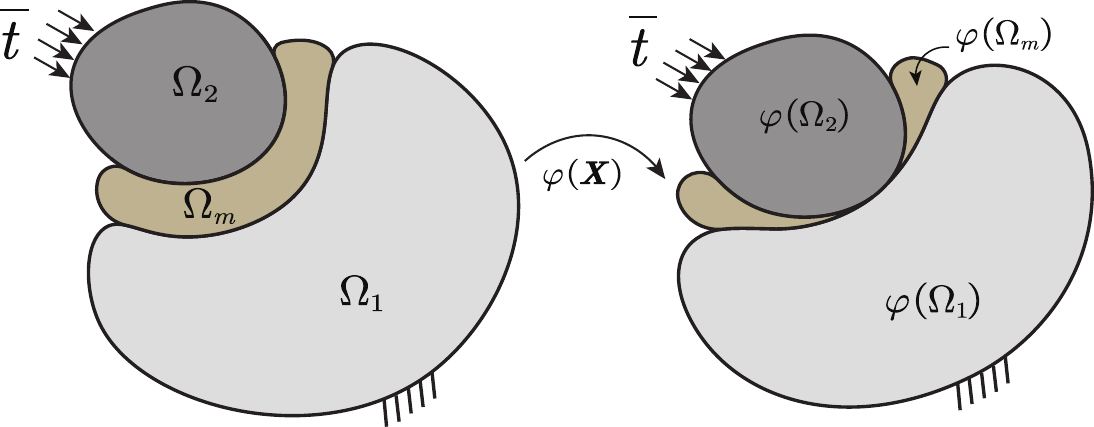}
\caption{Contact of two bodies with a third medium, $\phi$ is a map between the initial and current configuration.}
\label{p10.s2.tmc}
\end{figure}

The basic idea of the third medium contact is to introduce a third medium $\Omega_m$ between two bodies $\Omega_1$ and $\Omega_2$, as shown in Fig.\ref{p10.s2.tmc}.
Then, the contact problem can be formulated as a finite deformation problem.
The kinematic relation between the two bodies and the third medium are based on a formulation in the initial configuration.
Under the large deformation assumption,
the motion of the solid is then governed by the following equation
\begin{equation}
    \label{p10.s2.potential}
    \Pi(\bm{u}) = \int_\Omega\left[\Psi(\bm{F})-\bm{b}\cdot\bm{u}\right]\ud\Omega-\int_{\Gamma_N}\bm{\bar{t}}\cdot\bm{u}\ud\Gamma_N\Rightarrow STAT,
\end{equation}
where $\bm{u}$ is the displacement vector, $\bm{F}=\nabla\bm{u}+\bm{I}$ is the deformation gradient, 
$\Psi$ is the strain energy density function, $\bm{b}$ is the body force, and $\bm{\bar{t}}$ is the traction on the Neumann boundary $\Gamma_N$.

For compressible materials, a strain energy density function in Eq.\eqref{p10.s2.potential} can be expressed as
\begin{equation}
    \label{p10.s2.strainEnergy}
    \Psi(\bm{C}) = \frac{K}{2}\left(\ln J\right)^2+\frac{\mu}{2}\left(J^{-\frac{2}{3}}\text{tr}(\bm{C})-3\right) = \Psi_{vol}+\Psi_{iso},
\end{equation}
where $\bm{C}=\bm{F}^T\bm{F}$ is the right Cauchy-Green deformation tensor, $J=\det(\bm{F})$ is the volume change, 
$K$ is the bulk modulus, and $\mu$ is the shear modulus. 
For the bodies $\Omega_1$ and $\Omega_2$, other compressible or nearly incompressible strain energy densities can also be chosen.
Substuting the strain energy density function into the variation of Eq.\eqref{p10.s2.potential}, the weak form of the governing equations can be obtained.

\subsection{Strain energy density for the third medium contact}
\label{p10.s2.2}
The third medium should be a highly compressible material with the following strain energy density function
\begin{equation}
    \label{p10.s2.Psitmc}
    \Psi^{TMC}(\bm{u}) =\gamma \left[\Psi_m(\bm{u})+\alpha_r\Psi_r(\bm{u})\right] = \gamma\Psi_m(\bm{u})+\gamma\alpha_r\Psi_r(\bm{u}),
\end{equation}
where $\gamma$ and $\alpha_r$ are two parameters, $\Psi_m(\bm{u})$ is the strain energy density function, and $\Psi_r(\bm{u})$ is the regularization term to control the element distortions automatically.
By selecting a small parameter $\gamma$ and $\alpha_r$, the third medium can act as a highly soft material before contact 
and then rapidly provide a stiffness once the bodies getting close and contact.

In order to meet the above requirements, the strain energy density function $\Psi_m(\bm{u})$ for the third medium can be selected as 
\begin{equation}
    \label{p10.s2.Psim}
    \Psi_m(\bm{u}) = \frac{K}{2}\left(\ln J\right)^2+\frac{\mu}{2}\left(J^{-\frac{2}{3}}\text{tr}(\bm{C})-3\right).
\end{equation}
Since the contact between the solids is equivalent to the third medium being compressed to zero volume,
the term $\ln J$ has the property $J\rightarrow 0$ and $\Psi_m\rightarrow\infty$.
As discussed in \cite{TMC1, TMC2}, the second term in Eq.\eqref{p10.s2.Psim} for $J\rightarrow 0$  is sufficient to prevent penetration under plane strain conditions.
Then the strain energy density function for the third medium contact can be expressed as
\begin{equation}
	\Psi_m(\bm{u}) =\frac{\mu}{2}\left(J^{-\frac{2}{3}}\text{tr}(\bm{C})-3\right).
\end{equation}

\begin{figure}[htbp]
\centering
\includegraphics[width=1.0\textwidth]{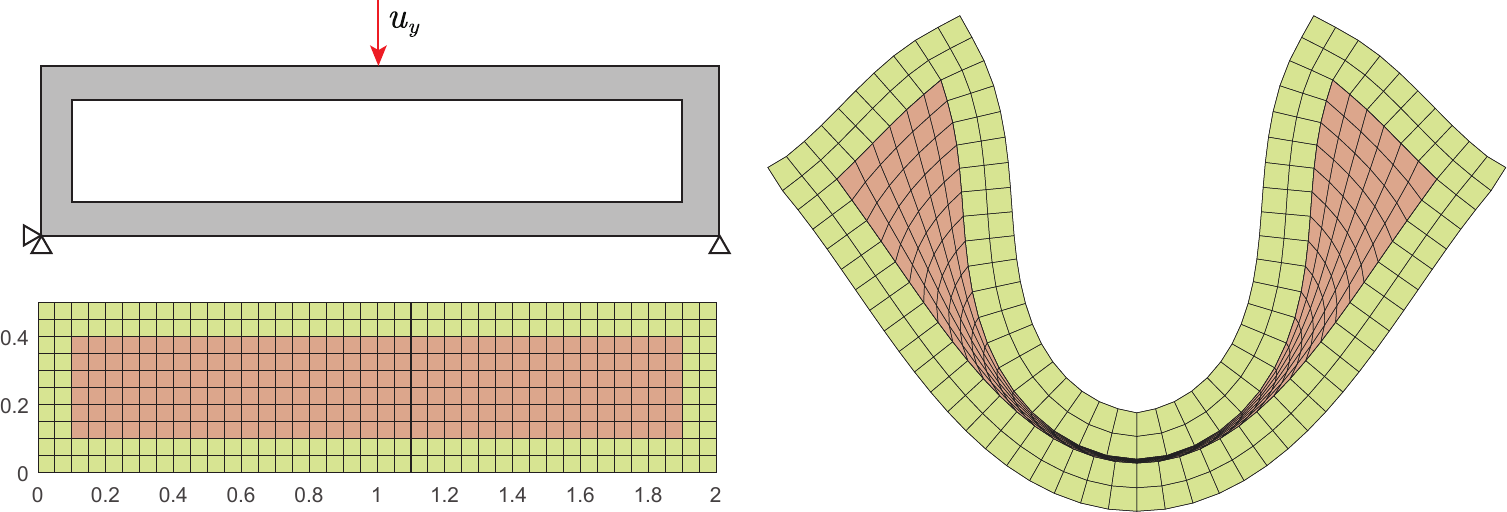}
\caption{Contact under finite strain assumption with third medium contact.}
\label{p10.s2.mesh}
\end{figure}

But it is also mentioned that, as shown in Fig.\ref{p10.s2.mesh}, 
the elements in the third medium is highly distorted when the two bodies are in contact (for detail, see Section \ref{p10.s5.1}).
Basically, the highly distorted elements will not allow to use the finite element method or virtual element method directly.

\subsection{Regularization term for the third medium contact}
\label{p10.s2.3}
In order to control the highly distorted elements in the third medium automatically, 
a regularization term $\Psi_r(\bm{u})$ is added as given in Eq.\eqref{p10.s2.Psim}.
As introduced in \cite{TMC3}, the regularization term can be expressed as
\begin{equation}
    \label{p10.s2.3.Psir0}
    \Psi_r(\bm{u}) = \frac{1}{2} \nabla \bm{F}\vdots\nabla\bm{F} = \frac{1}{2} \mathbb{H}\bm{u}\vdots \mathbb{H}\bm{u},
\end{equation}
where  $\nabla\bm{F}$ is the gradient of the deformation gradient,
\begin{equation}
	\nabla\bm{F}\vdots \nabla\bm{F} = \frac{\partial^2 u_i}{\partial X_j\partial X_k}\frac{\partial^2 u_i}{\partial X_j\partial X_k}.
\end{equation}
Besides, the regularization term can be multiplied with the scaling term $e^{-\beta |F|}$ as 
\begin{equation}
    \label{p10.s2.3.Psir}
    \Psi_r(\bm{u}) = \frac{1}{2} e^{-\beta |F|}\nabla \bm{F}\vdots\nabla\bm{F}.
\end{equation}

As mentioned in \cite{TMC1}, the regularization term $\Psi_r(\bm{u})$ cannot accurately simulate certain specific contacts and cannot converge when the mesh is deformed.
Then,  there are some other regularization terms proposed in the literature \cite{TMC1, TMC2, TMC3}.
Faltus et al.\cite{TMC1} proposed a new regularization term by subtracting a new term from the Hu-Hu regularization (Eq.\eqref{p10.s2.3.Psir0}) 
\begin{equation}
    \label{p10.s2.3.Psir1}
    \Psi_r = \frac{1}{2}\left[\nabla\bm{F}\vdots\nabla\bm{F} - \frac{1}{n_{dim}} \text{Div}(\nabla\bm{u})\cdot\text{Div}(\nabla\bm{u})\right],
\end{equation}
where $n_{dim}$ is the number of dimensions, $\text{Div}(\nabla\bm{u}) = \frac{\partial^2 u_i}{\partial X_j\partial X_j}$ is the divergence of the displacement gradient tensor.

Besides, the regularization term can also be expressed in terms of the gradient of rotation tensor $\nabla\bm{R}$ as \cite{TMC1}
\begin{equation}
    \label{p10.s2.3.Psir2}
    \Psi_r = \frac{1}{2}\left(\nabla\bm{R}\vdots\nabla\bm{R} + \nabla J\cdot\nabla J\right),
\end{equation}
where $\bm{R}$ is the rotation tensor and $\nabla J\cdot\nabla J = \frac{\partial J}{\partial X_j}\frac{\partial J}{\partial X_j}$.
As introduced in Wriggers et al.\cite{TMC2}, the rotation tensor $\bm{R}$ can be expressed as
\begin{equation}
	\bm{R} = \begin{bmatrix}
		\cos\varphi & -\sin\varphi & 0\\
		\sin\varphi & \cos\varphi & 0\\
        0 & 0 & 1
	\end{bmatrix}
\end{equation}
for 2D problems, where $\varphi$ is the rotation angle. 
Then, the rotation angle $\varphi$ can be solved by considering the symmetry of the right stretch tensor $\bm{U} = \bm{R}^T\bm{F}$
\begin{equation}
    U_{12} = U_{21} \Rightarrow R_{11}F_{12}+R_{21}F_{22} = R_{12}F_{11}+R_{22}F_{21} \Rightarrow \tan \varphi = \left[\frac{F_{12}-F_{21}}{F_{11}+F_{22}}\right],
\end{equation}
and $\varphi = \arctan(\frac{F_{12}-F_{21}}{F_{11}+F_{22}})$.
Furthermore, the gradient of the rotation tensor can be expressed as
\begin{equation}
	\nabla\bm{R} = \frac{\partial\bm{R}}{\partial\varphi}\otimes\nabla\varphi
\end{equation}
which lead to 
\begin{equation}
    \label{p10.s2.3.nablaR}
	\nabla\bm{R}\vdots\nabla\bm{R} = \frac{\partial\bm{R}}{\partial\varphi}\otimes\nabla\varphi:\frac{\partial\bm{R}}{\partial\varphi}\otimes\nabla\varphi
	= \left(\frac{\partial\bm{R}}{\partial\varphi}\vdots\frac{\partial\bm{R}}{\partial\varphi}\right)\nabla\varphi\cdot\nabla\varphi
	=2 \nabla\varphi\cdot\nabla\varphi.
\end{equation}

Substuting Eq.\eqref{p10.s2.3.nablaR} into the regularization term Eq.\eqref{p10.s2.3.Psir2}, the regularization term can be expressed as
\begin{equation}
    \label{p10.s2.3.rw1}
    \Psi_r = \frac{1}{2}\left(\nabla\varphi\cdot\nabla\varphi+\nabla J\cdot\nabla J\right).
\end{equation}
Considering $\tan \varphi = \left[\frac{F_{12}-F_{21}}{F_{11}+F_{22}}\right]$, another choice for the regularization term can be selected as
\begin{equation}
    \label{p10.s2.3.rw2}
	\Psi_r = \frac{1}{2}
	\left[
		\nabla\left[\frac{F_{12}-F_{21}}{F_{11}+F_{22}}\right]\cdot\nabla\left[\frac{F_{12}-F_{21}}{F_{11}+F_{22}}\right]+\nabla J\cdot\nabla J
	\right].
\end{equation}
The gradient of the Jacobian $\nabla J$ in Eqs.\eqref{p10.s2.3.rw1} and \eqref{p10.s2.3.rw2}  can be calculated by
\begin{equation}
	\nabla J = \frac{\partial\det\bm{F}}{\partial\bm{F}}:\nabla\bm{F} = \det\bm{F}\cdot\bm{F}^{-T}:\nabla\bm{F}.
\end{equation}
Different regularization terms are summarized in Table \ref{p10.s2.3.t1}.
Besides, the scaling term $e^{-\beta |F|}$ in Eq.\eqref{p10.s2.3.Psir} can also be used for the above regularization terms.

\begin{table}[htbp]
    \centering
    \caption{Different regularization terms for the third medium contact, $n_{dim}$ is the dimension.}
    \label{p10.s2.3.t1}
    \begin{tabular}{cc}
    \toprule          
    Reference       & $\Psi_r$ \\
    \midrule
    \cite{TMC1} & $\frac{1}{2} \nabla \bm{F}\vdots\nabla\bm{F}$ \\
    \cite{TMC3} & $\frac{1}{2} \left[\nabla\bm{F}\vdots\nabla\bm{F} - \frac{1}{n_{dim}} \text{Div}(\nabla\bm{u})\cdot\text{Div}(\nabla\bm{u})\right]$ \\
    \cite{TMC2} & $\frac{1}{2} \left(\nabla\varphi\cdot\nabla\varphi+\nabla J\cdot\nabla J\right)$ \\
    \cite{TMC2} & $\frac{1}{2} \left[\nabla\left[\frac{F_{12}-F_{21}}{F_{11}+F_{22}}\right]\cdot\nabla\left[\frac{F_{12}-F_{21}}{F_{11}+F_{22}}\right]+\nabla J\cdot\nabla J\right]$\\
    \bottomrule  
    \end{tabular}
\end{table}

Different regularization terms for the third medium contact will be used and comprared in the numerical examples.
Substuting the strain energy density function $\Psi_m(\bm{u})$ and the regularization term $\Psi_r(\bm{u})$ into Eq.\eqref{p10.s2.potential}, the potential energy for the third medium contact can be expressed as
\begin{equation}
    W^{TMC}(\bm{u}) = \int_{\Omega_m}\gamma \left[\Psi_c(\bm{u})+\alpha_r\Psi_r(\bm{u})\right]\ud\Omega.
\end{equation}
Considering the variation of the total energy, a nonlinear system of equation can be obtained and a linearization is necessary.
In addition, since higher-order derivatives of the displacement are taken into account, a smaller time step is required for the calculation.

\section{Linearization and tangent stiffness matrix}
The third medium contact theory is based on the framework of hyperelasticity, 
so it's necessary to consider the variation of potential energy and its linearization process.
Easy to find that the potential energy for the hyperelastic bodies $\Omega_1$ and $\Omega_2$ is a functional of $\bm{C}$ (or $\bm{F}$),
and the potential energy for third medium is a functional of $\bm{F}$ and $\nabla \bm{F}$. 
This difference leads to different treatments for different domains.
\label{p10.s3}
\subsection{Linearization for the hyperelastic body}
\label{p10.s3.1}
For the hyperelastic bodies $\Omega_1$ and $\Omega_2$, the strain energy density function in Eq.\eqref{p10.s2.strainEnergy} can be selected and the potential energy has the form as 
\begin{equation}
    \label{p10.s3.1.W}
    W(\bm{C}) = \int_{\Omega}\Psi(\bm{C})\ud\Omega.
\end{equation}
The variational form of the principle of the strain energy yields
\begin{equation}
	\delta W = \int_{\Omega}	\frac{\partial \Psi}{\partial\bm{E}}: \delta\bm{E}\ud\Omega = \int_{\Omega}\bm{S}:\delta\bm{E}\ud\Omega,
\end{equation}
where $\bm{S}$ is the second Piola-Kirchhoff stress and $\bm{E}$ is the Green-Lagrange strain tensor
\begin{equation}
    \bm{E} = \frac{1}{2}\left(\bm{F}^T\cdot\bm{F}-\bm{I}\right),
\end{equation}
where $\bm{I}$ is the second order identity tensor.
Then, the increasement can be expressed as
\begin{equation}
    \label{p10.s3.1.ddw}
	\begin{aligned}
        \Delta\delta W &= \int_{\Omega}\delta\bm{E}:\frac{\partial^2\Psi}{\partial\bm{E}\partial\bm{E}}:\Delta\bm{E}\ud\Omega+\int_\Omega\Delta(\delta\bm{E}):\bm{S}\ud\Omega\\
        & = \int_{\Omega}\delta\bm{E}:\mathbb{D}:\Delta\bm{E}\ud\Omega+\int_\Omega\Delta(\delta\bm{E}):\bm{S}\ud\Omega
    \end{aligned}
\end{equation}
where $\mathbb{D}$ is the constitutive tensor. 

Considering the strain energy density function in Eq.\eqref{p10.s2.strainEnergy},  the second Piola-Kirchhoff stress can be expressed as
\begin{equation}
	\label{p10.s3.1.S}
	\begin{aligned}
		\bm{S} =2K\ln J\frac{1}{J}\frac{\partial J}{\partial\bm{C}}
	-\frac{2}{3}\mu J^{-5/3}\frac{\partial J}{\partial\bm{C}}\text{tr}(\bm{C})
	+\mu J^{-2/3}\bm{I}.\\
	\end{aligned}
\end{equation}
Considering $J = \sqrt{\text{det}\bm{C}}$,  the second Piola-Kirchhoff stress has the last form as
\begin{equation}
    \label{p10.s3.1.S2}
	\bm{S} = K\ln J\bm{C}^{-1}-\frac{1}{3}\mu J^{-2/3}\text{tr}(\bm{C})\bm{C}^{-1}+\mu J^{-2/3}\bm{I}.
\end{equation}

Besides, the constitutive tensor $\mathbb{D} = 2\partial\bm{S}/\partial\bm{C}$ has the form as
\begin{equation}
	\label{p10.s3.1.D}
	\begin{aligned}
		\mathbb{D}_{ijkl} =& K\left(\bm{C}^{-1}\right)_{ij}\left(\bm{C}^{-1}\right)_{kl} - 2K\ln J\mathcal{I}_{ijkl}\\
		&+ \frac{2\mu}{3}J^{-2/3}\left(
			-\delta_{ij}\left(\bm{C}^{-1}\right)_{kl}-\delta_{kl}\left(\bm{C}^{-1}\right)_{ij}
			+\frac{1}{3}\text{tr}(\bm{C})\left(\bm{C}^{-1}\right)_{ij} \left(\bm{C}^{-1}\right)_{kl}
			+\text{tr}(\bm{C})\mathcal{I}_{ijkl}
		\right)
	\end{aligned}
\end{equation}
with the matrix expressed as
\begin{equation}
	\begin{aligned}
		\mathbb{D} =& K\bm{C}^{-1}\otimes\bm{C}^{-1} - 2K\ln J\bm{\mathcal{I}}\\
		& +\frac{2\mu}{3}J^{-2/3}\left(
		-\bm{I}\otimes\bm{C}^{-1}-\bm{C}^{-1}\otimes\bm{I}
		+\frac{1}{3}\text{tr}(\bm{C})\bm{C}^{-1}\otimes\bm{C}^{-1}
		+\text{tr}(\bm{C})\bm{\mathcal{I}}
	\right),
	\end{aligned}
\end{equation}
where 
\begin{equation}
	\mathcal{I}_{ijkl} = -\frac{\partial\left(\bm{C}^{-1}\right)_{ij}}{\partial C_{kl}} 
	= \frac{1}{2}\left[\left(\bm{C}^{-1}\right)_{ik}\left(\bm{C}^{-1}\right)_{jl} +\left(\bm{C}^{-1}\right)_{il}\left(\bm{C}^{-1}\right)_{jk}\right].
\end{equation}

The detailed derivation process can be found in \cite{Wriggers2008}.
In order to simplify the above complex derivation process,
the second Piola-Kirchhoff stress and the constitutive tensor can also be obtained by the automatic differentiation, see \cite{Korelc2016}.
The automatic differentiation matlab code for the second Piola-Kirchhoff stress and the constitutive tensor is given in the appendix \ref{p10.appendix1}.

\subsection{Linearization for the third medium}
\label{p10.s3.2}
As discussed before, the potential energy for the third medium is a functional of $\bm{F}$ and $\nabla \bm{F}$ as
\begin{equation}
	W^{TMC}(\bm{F},\nabla\bm{F}) = \int_{\Omega_m}\Psi^{TMC}(\bm{F},\nabla\bm{F})\ud\Omega.
\end{equation}
Then, the variation of the potential energy should be
\begin{equation}
	\begin{aligned}
        \delta W^{TMC} &= \int_{\Omega_m}
	\left[\frac{\partial \Psi^{TMC}}{\partial\bm{F}}:\delta\bm{F}+\frac{\partial \Psi^{TMC}}{\partial\nabla\bm{F}}\vdots\delta\nabla\bm{F}\right]
	\ud\Omega\\
    & = \int_{\Omega_m}
	\left(\bm{P}:\delta\bm{F}+\bm{T}\vdots\delta\nabla\bm{F}\right)
	\ud\Omega,
    \end{aligned}
\end{equation}
where 
\begin{equation}
    \bm{P} = \frac{\partial \Psi^{TMC}}{\partial\bm{F}},\quad 
    \bm{T} = \frac{\partial \Psi^{TMC}}{\partial\nabla\bm{F}}.
\end{equation}

For linearization, we consider the second variation of the potential energy
\begin{equation}
    \label{p10.s3.2.W2}
	\begin{aligned}
		\delta^2 W^{TMC} &=\int_{\Omega_m}
		\delta\bm{F}:\frac{\partial^2\Psi^{TMC}}{\partial\bm{F}\partial\bm{F}}:\delta\bm{F}+
		\delta\nabla\bm{F}\vdots\frac{\partial^2\Psi^{TMC}}{\partial\nabla\bm{F}\partial\bm{F}}:\delta\bm{F}\\
		&
		+\delta\bm{F}:\frac{\partial^2\Psi^{TMC}}{\partial\bm{F}\partial\nabla\bm{F}}\vdots\delta\nabla\bm{F}
		+\delta\nabla \bm{F}\vdots\frac{\partial^2\Psi^{TMC}}{\partial\nabla\bm{F}\partial\nabla\bm{F}}\vdots\delta\nabla\bm{F}
		\ud\Omega\\
        &=\int_{\Omega_m}\left(
            \delta\bm{F}:\mathbb{D}:\delta\bm{F}
        +\delta\nabla\bm{F}\vdots\mathbb{A}:\delta\bm{F}
        +\delta\bm{F}:\mathbb{A}\vdots\delta\nabla\bm{F}
        +\delta\nabla \bm{F}\vdots\mathbb{B}\vdots\delta\nabla\bm{F}
        \right)\ud\Omega,
	\end{aligned}
\end{equation}
where 
\begin{equation}
	\mathbb{D} = \frac{\partial^2\Psi^{TMC}}{\partial\bm{F}\partial\bm{F}},\quad
	\mathbb{A} = \frac{\partial^2\Psi^{TMC}}{\partial\nabla\bm{F}\partial\bm{F}},\quad
	\mathbb{B} = \frac{\partial^2\Psi^{TMC}}{\partial\nabla\bm{F}\partial\nabla\bm{F}}.
\end{equation}

For different regularization term as summaried in Table \ref{p10.s2.3.t1}, the derivation of specific explicit expressions for tensors $\bm{P},\bm{T}$ and $\mathbb{A},\mathbb{B},\mathbb{D}$ would be very complicated.
For example, for the regularization term mentioned in Eq.\eqref{p10.s2.3.Psir0},$\Psi_r = \frac{1}{2}\nabla\bm{F}\vdots\nabla\bm{F}$, revelent tensors have the form as 
\begin{equation}
    \label{p10.s3.2.PT}
    \bm{P} = \frac{\partial \Psi_r}{\partial\bm{F}} = 0,\quad 
    \bm{T} = \frac{\partial \Psi_r}{\partial\nabla\bm{F}} = \nabla\bm{F},
\end{equation}
\begin{equation}
    \label{p10.s3.2.DAB}
    \mathbb{D} = 0,\quad \mathbb{A} = 0,\quad 
    \mathbb{B} = \frac{\partial^2\Psi_r}{\partial\nabla\bm{F}\partial\nabla\bm{F}} = \mathbb{I},
\end{equation}
where $\mathbb{I}$ is a six-order tensor $\mathbb{I} = \delta_{im}\delta_{ij}\delta_{kp}$.
Substuting Eqs.\eqref{p10.s3.2.PT} and \eqref{p10.s3.2.DAB} into Eq.\eqref{p10.s3.2.W2}, the second variation of the potential energy can be obtained lastly.
But unfortunately, if the scaling term $e^{-\beta |F|}$ is considered or other regularization terms are selected, it could be hard to get the specific explicit expressions of above tensors.

To simplify, the tensors $\bm{P},\bm{T}$ and $\mathbb{A},\mathbb{B},\mathbb{D}$ can be obtained by the automatic differentiation in MATLAB or the software tool AceGen \cite{Korelc2016}.
For example, for the regularization term discussed in Eq.\eqref{p10.s2.3.Psir1} with the scaling term
\begin{equation}
    \Psi_r = \frac{1}{2} e^{-\beta |F|}\left[\nabla\bm{F}\vdots\nabla\bm{F} - \frac{1}{n_{dim}} \text{Div}(\nabla\bm{u})\cdot\text{Div}(\nabla\bm{u})\right],
\end{equation}
the MATLAB code for the tensors $\bm{P},\bm{T}$ and $\mathbb{A},\mathbb{B},\mathbb{D}$ is given in the appendix \ref{p10.appendix2}.

For the hyperelastic bodies $\Omega_1$ and $\Omega_2$, the above framework and code can also be used with $\mathbb{A} = 0,\mathbb{B} = 0$ and $\bm{T} = 0$.
Then the first Piola-Kirchhoff stress $\bm{P}$ is used rather than the second Piola-Kirchhoff stress $\bm{S}$ as gievn in Eq.\eqref{p10.s3.1.S2}.

\section{Stabilization-free virtual element method for third medium contact}
The third medium contact is based on the framework of hyperelasticity and the regularization term contains higher-order derivatives of displacement.
Then, the finite element method or virtual element method with higher-order ansatz is necessary.
In order to avoid constructing the projection of the second-order gradient, 
we can construct a virtual element method without unstable terms to indirectly calculate the regularization term.

\label{p10.s4}
\subsection{High-order gradient projection operator}
\label{p10.s4.1}
In the third medium contatc, the gradient of the deformation gradient $\nabla\bm{F}$ is needed so the second-order virtual element method is necessary.
We introduce a projection operator $\Pi_{k,E}^\nabla:\mathcal{V}_k(E)\rightarrow \mathbb{P}_k(E), v\mapsto \Pi_{k,E}^\nabla v$ by
\begin{equation}
    \label{p07.s3.1.energy}
    \left\{
        \begin{aligned}
            &\int_E\nabla\left(\Pi_{k,E}^\nabla v-v\right)\cdot\nabla p_k\ud\Omega = 0,\quad \forall v\in H^1(E), p_k\in\mathbb{P}_k(E),\\
            &\frac{1}{|E|} \int_E \left(\Pi_{k,E}^\nabla v-v\right)\ud\Omega = 0,
        \end{aligned}
    \right.
\end{equation}
where $\mathcal{V}_k(E)$ is the local virtual element space
\begin{equation}
	\mathcal{V}_k(E) :=\left\{v\in H^1(E):\Delta u\in\mathbb{P}_{k-2}(E)\quad\text{in}\quad E,\quad u|_{\partial E} = \mathcal{B}_k(\partial E)\right\},
\end{equation}
\begin{equation}
	\mathcal{B}_k(\partial E):=\left\{v\in C(\partial E):v_e\in\mathbb{P}_k(e),\quad e\subset \partial E\right\}.
\end{equation}

In $\mathcal{V}_k(E)$, the degrees of freedom are selected as 
\begin{itemize}
	\item $\bm{\chi}^1(E)$: for $k\geq 1$,the values of $v$ at the vertices;
	\item $\bm{\chi}^2(E)$: for $k>1$, the values of $v$ at $k-1$ uniformly spaced points on each edge $e$;
	\item $\bm{\chi}^3(E)$: for $k>1$, the moments 
	\begin{displaymath}
		\frac{1}{|E|}\int_E v p_{k-2}\ud\Omega,\quad\forall p_{k-2}\in \mathbb{P}_{k-2}(E).
	\end{displaymath}
\end{itemize}
Given the above degrees of freedom, the projection operator $\Pi_{k,E}^\nabla$ can be solved based on Eq.\eqref{p07.s3.1.energy}. 
For detail, see \cite{Veiga2014}.
When using this projection operator to discretize variables, 
we need to construct additional stabilization terms and additional projection operators to calculate higher-order derivatives.

Inspired by Ref.\cite{StabVEM6}, a local enlarged enhancement virtual element space of order $l$ is introduced 
based on the higher order $L_2$ polynomial projection 
\begin{equation}
    \Pi_{l,E}^0\nabla:H^1(E)\rightarrow\left[\mathbb{P}_l(E)\right]^2.
\end{equation}
The selection of $l$ is discussed in \cite{StabVEM1, StabVEM6} for the first-order stabilization-free VEM.
For the higher order method $k>1$, by considering the element eigenvalue problem for plane elasticity \cite{StabVEM3}, a sufficient inequality is given by 
\begin{equation}
    N_E\leq 2l-2k+5,
\end{equation}
where $N_E$ is the number of vertices of element $E$.

To solve the $L_2$ projection operator, the orthogonality condition should be satisfied as
\begin{equation}
    \label{p10.s4.1.L2}
    \int_E\bm{p}\cdot\Pi_{l,E}^0\nabla v\ud\Omega = \int_E\bm{p}\cdot\nabla v\ud\Omega,\quad \bm{p}\in \left[\mathbb{P}_l(E)\right]^2.
\end{equation}
Expanding the right side of Eq.\eqref{p10.s4.1.L2} yields
\begin{equation}
    \label{p10.s4.1.L2.2}
    \int_E\bm{p}\cdot\nabla v\ud\Omega=\int_{\partial E}\bm{p}\cdot\bm{n}_e v\ud\Gamma-\int_E\left(\text{div}\bm{p}\right)v\ud\Omega,
\end{equation}
where $\bm{n}_e$ is the out normal of edge $e$.
The last term in \eqref{p10.s4.1.L2.2} is computable 
as $\int_E\left(\text{div}\bm{p}\right)v\ud\Omega=\int_E\left(\text{div}\bm{p}\right)\Pi_{1,E}^\nabla v\ud\Omega$.
Then the $L_2$ projection operator $\Pi_{l,E}^0\nabla$ can be solved based on the degrees of freedom used for $\mathcal{V}_k(E)$.

By solving the $L_2$ projection operator $\Pi_{l,E}^0\nabla:H^1(E)\rightarrow\left[\mathbb{P}_l(E)\right]^2$,
the gradient of the variable $u_h$ can be approximated as follows
\begin{equation}
    \label{p10.s4.2.L2}
    \nabla v= \Pi_{l,E}^0\nabla v =  \left(\bm{N}^p\right)^T\bm{\Pi}^m\tilde{\bm{v}},
\end{equation}
where $\tilde{\bm{v}}$ is a vector of node variable value in $E$,
$\bm{\Pi}^m$ is the matrix representation of the $L_2$ projection operator $\Pi_{l,E}^0\nabla$.
Besides, $\bm{N}^p$ in Eq.\eqref{p10.s4.2.L2} is
\begin{equation}
    \left(\bm{N}^p\right)^T :=
    \begin{bmatrix}
        1 & \xi & \eta& \cdots& \eta^l & 0&0&0&\cdots &0 \\
        0&0&0&\cdots &0& 1 & \xi & \eta& \cdots& \eta^l
    \end{bmatrix} = 
	\begin{bmatrix}
		\bm{m}_l^T & \\ & \bm{m}_l^T
	\end{bmatrix},
\end{equation}
and $\bm{m}_l$ is the basic for the scaled polynomials of order $l$. 

Then the second-order derivative of the variable $v$ can also be calculated by 
\begin{equation}
    \label{p10.s4.1.ddv}
    \nabla\nabla v = \left[\nabla \left(\bm{N}^p\right)^T\right]\bm{\Pi}^m\tilde{\bm{v}}.
\end{equation}
Compared with FEM, VEM does not require coordinate transformation, so the calculation of the second-order derivative is simpler. 
If FEM is used, the calculation of the second-order derivatives of the variables will be more complicated.

\subsection{Discrete formualtion}
\label{p10.s4.2}
We consider the hyperelastic bodies $\Omega_1,\Omega_2$ and the third medium $\Omega_m$ as a unified body, so the strain energy density function $\Psi(\bm{F},\nabla\bm{F})$ is considered.
As discussed before, the potential energy is 
\begin{equation}
	W(\bm{F},\nabla\bm{F}) = \int_\Omega\Psi(\bm{F},\nabla\bm{F})\ud\Omega.
\end{equation}
The variation should be
\begin{equation}
	\delta W = \int_{\Omega_m}
	\left[\frac{\partial \Psi}{\partial\bm{F}}:\delta\bm{F}+\frac{\partial \Psi}{\partial\nabla\bm{F}}\vdots\delta\nabla\bm{F}\right]
	\ud\Omega,
\end{equation}
and the second variation of the potential energy
\begin{equation}
    \delta^2 W = \int_\Omega\left(
            \delta\bm{F}:\mathbb{D}:\delta\bm{F}
        +\delta\nabla\bm{F}\vdots\mathbb{A}:\delta\bm{F}
        +\delta\bm{F}:\mathbb{A}\vdots\delta\nabla\bm{F}
        +\delta\nabla \bm{F}\vdots\mathbb{B}\vdots\delta\nabla\bm{F}
        \right)\ud\Omega.
\end{equation}

In the framework of third medium contact, it is necessary to calculate the deformation gradient $\bm{F}$ and the gradient of the deformation gradient $\nabla\bm{F}$.
To facilitate matrix operations, the variation of the deformation gradient $\bm{F}$ can be arranged into a first-order vector as 
\begin{equation}
	\begin{aligned}
		\delta\hat{\bm{F}} = \begin{bmatrix}
			\frac{\partial\delta u_x}{\partial x} & \frac{\partial\delta u_y}{\partial x} & 
			\frac{\partial\delta u_x}{\partial y} & \frac{\partial\delta u_y}{\partial y}
		\end{bmatrix}^T
		= \underbrace{\bm{A}\left(\bm{N}_p\right)^T\bm{\Pi}_m}_{\bm{B}_1}\delta\tilde{\bm{u}},
	\end{aligned}
\end{equation}
where
\begin{equation}
        \bm{A} = \begin{bmatrix}
        1 & 0 & 0 & 0\\
        0 & 0 & 0 & 1\\
        0 & 1 & 1 & 0\\
    \end{bmatrix}, \bm{N}_p = \bm{N}^p\otimes\mathbb{I}_2,\bm{\Pi}_m = \bm{\Pi}^m\otimes\mathbb{I}_2,
    \bm{B}_1 = \bm{A}\left(\bm{N}_p\right)^T\bm{\Pi}_m,
\end{equation}
and $\tilde{\bm{u}} = [\tilde{\bm{u}}_x,\tilde{\bm{u}}_y]^T$ is the node displacement vector.
Besides, $\otimes$ is the Kronecker product and $\mathbb{I}_2$ represents $2\times 2$ order identity matrix.

Besides, the variation of the gradient of the deformation gradient $\nabla\bm{F}$ can be arranged into a first-order vector by column
\begin{equation}
	\widehat{\nabla\delta\bm{F}}  = 
	\begin{bmatrix}
		\frac{\partial^2\delta u_x}{\partial x\partial x} & 
		\frac{\partial^2\delta u_y}{\partial x\partial x} &
		\frac{\partial^2\delta u_x}{\partial y\partial x} & 
		\frac{\partial^2\delta u_y}{\partial y\partial x} & 
		\frac{\partial^2\delta u_x}{\partial x\partial y} & 
		\frac{\partial^2\delta u_y}{\partial x\partial y} & 
		\frac{\partial^2\delta u_x}{\partial x\partial x} & 
		\frac{\partial^2\delta u_y}{\partial x\partial x}
	\end{bmatrix}^T = \bm{B}_2\delta\tilde{\bm{u}},
\end{equation}
where the matrix $\bm{B}_2$ can be computed based on Eq.\eqref{p10.s4.1.ddv}.

Based on the above definition, 
the matrix formualtion of the variation of the potential energy in the stabilization-free virtual element method is 
\begin{equation}
	\delta W = \delta\tilde{\bm{u}}^T\int_\Omega \left(\bm{B}_1^T\hat{\bm{P}}+\bm{B}_2^T\hat{\bm{T}}\right) \ud\Omega,
\end{equation} 
where 
\begin{equation}
	\hat{\bm{P}} = \begin{bmatrix}
		P_{11} & P_{21} & P_{12} & P_{22}
	\end{bmatrix}^T,
\end{equation}
\begin{equation}
	\hat{\bm{T}} = \begin{bmatrix}
		T_{111} & T_{211} & T_{121} & T_{221} & T_{112} & T_{212} & T_{122} & T_{222}
	\end{bmatrix}^T.
\end{equation}

Besides, the second variation of the potential energy is 
\begin{equation}
    \label{p10.s4.2.ddw}
	\delta^2 W = \delta\tilde{\bm{u}}^T\left[
		\int_\Omega
	\left(
		\bm{B}_1^T\hat{\mathbb{D}}\bm{B}_1+\bm{B}_2^T\hat{\mathbb{A}}\bm{B}_1+\bm{B}_1^T\hat{\mathbb{A}}^T\bm{B}_2+\bm{B}_2^T\hat{\mathbb{B}}\bm{B}_2
	\right)\ud\Omega
	\right]\delta\tilde{\bm{u}},
\end{equation}
where 
the matrices $\hat{\mathbb{D}},\hat{\mathbb{A}},\hat{\mathbb{B}}$ are the matrix forms of the fourth-order tensor $\mathbb{D}$, the sixth-order tensor $\mathbb{A}$ and the eighth-order tensor $\mathbb{B}$,
and are given in the appendix \ref{p10.appendix3}.

Lastly, the tangent stiffness matrix on element level can be obtained based on the SFVEM as 
\begin{equation}
    \label{p10.s4.2.Ke}
    \bm{K}_E = \int_\Omega
	\left(
		\bm{B}_1^T\hat{\mathbb{D}}\bm{B}_1+\bm{B}_2^T\hat{\mathbb{A}}\bm{B}_1+\bm{B}_1^T\hat{\mathbb{A}}^T\bm{B}_2+\bm{B}_2^T\hat{\mathbb{B}}\bm{B}_2
	\right)\ud\Omega.
\end{equation}

For the hyperelastic bodies $\Omega_1$ and $\Omega_2$, the tangent stiffness matrix can be simplified as 
\begin{equation}
    \bm{K}_E = \int_\Omega
	\bm{B}_1^T\hat{\mathbb{D}}\bm{B}_1\ud\Omega.
\end{equation}
Besides, the tangent stiffness matrix for $\Omega_1$ and $\Omega_2$ can also obtained based on Eq.\eqref{p10.s3.1.ddw} (the second Piola-Kirchhoff stress is used), see Ref.\cite{StabXu1} for detail.

\section{Numerical examples}
\label{p10.s5}
In this section, we will present some numerical examples to verify the effectiveness of the SFVEM for the third medium contact.

\subsection{Self-contact within a box}
\label{p10.s5.1}
\begin{figure}[htbp]
    \centering
    \includegraphics[width=1.0\textwidth]{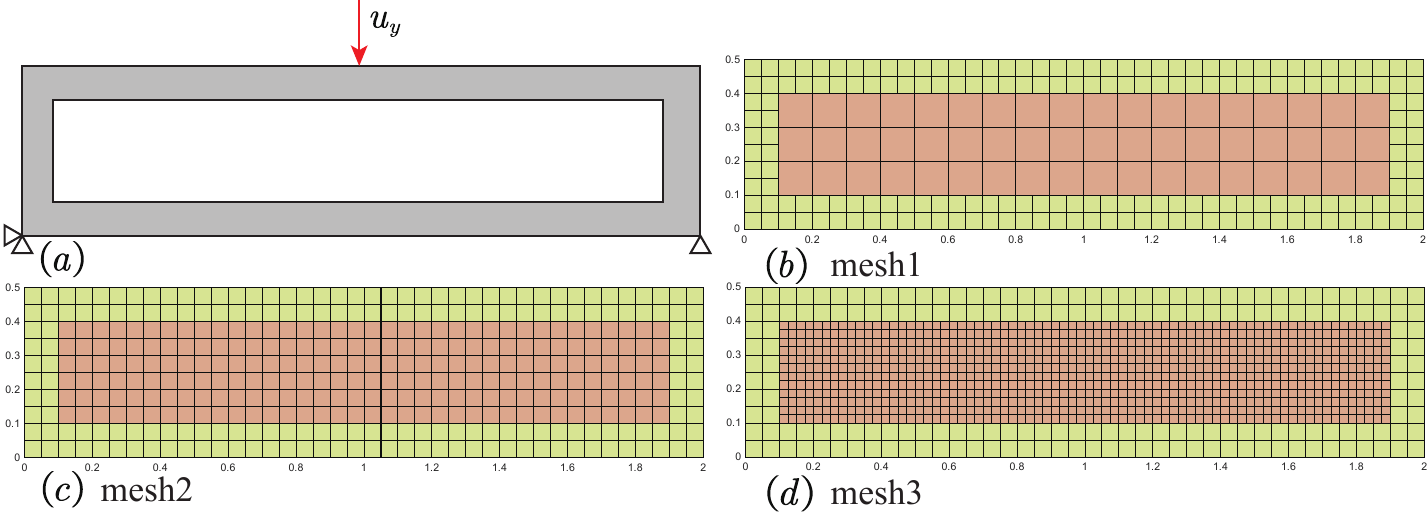}
    \caption{Self-contact within a box, (a) geometry model and boundary condition, (b)-(d) different meshes for box with third medium.}
    \label{p10.s5.1.f1}
\end{figure}

We consider a self-contact problem of a box with the third medium contact.
As shown in Fig.\ref{p10.s5.1.f1} (a), the box has a length of $L=2$ and a height of $H=0.5$.
Besides, a constant thickness of $T=0.1$ is set for the box.
The box is fixed at the lower left corner and simply supported at the lower right corner.
A vertical displacement load $u_y$ is applied in the middle at the top of the upper box.
The strain energy density function is given in Eq.\eqref{p10.s2.strainEnergy} with bulk modulus $K=20$ and shear modulus $\mu=10$.

Since the SFVEM allows the use of polygonal elements, it offers a significant advantage when handling meshes with hanging nodes.
This enables us to employ meshes of different densities in the third medium without requiring any special treatment in the other domains.
In this example, we consider three different mesh configurations for the third medium, as illustrated in Fig.~\ref{p10.s5.1.f1}(b)-(d).
The second-order SFVEM is employed and the order of the $L_2$ projection operator is set as $l=3$.
The vertical displacement load is selected as $u_y = -1.0$ with $100$ loading steps within the Newton-Raphson solution scheme.

Firstly, the performance of different parameters for the third medium is investigated.
Here the regularization term is selected as
$\Psi_r = \frac{1}{2} e^{-\beta |F|}\left(\nabla\bm{F}\vdots\nabla\bm{F} - \frac{1}{n_{dim}} \text{Div}(\nabla\bm{u})\cdot\text{Div}(\nabla\bm{u})\right)$ as mentioned in Eq.\eqref{p10.s2.3.Psir1},
where $\beta$ is selected as $\beta = 5$.
Then the potential energy for the third medium domain can be expressed as
\begin{equation}
    W^{TMC}(\bm{u}) = \int_{\Omega_m}\gamma \left[\frac{\mu}{2}\left(J^{-\frac{2}{3}}\text{tr}(\bm{C})-3\right)
    +\alpha_r \frac{1}{2} e^{-5 |F|}\left(\nabla\bm{F}\vdots\nabla\bm{F} - \frac{1}{2} \text{Div}(\nabla\bm{u})\cdot\text{Div}(\nabla\bm{u})\right)\right]\ud\Omega,
\end{equation}
and there are two parameters $\gamma$ and $\alpha_r$ in the potential energy of the third medium.
For different parameters $\gamma$ and $\alpha_r$, 
the gap $g$ between the surfaces of the upper and lower flange is recorded in Table \ref{p10.s5.1.t1}.
Easy to find that the gap $g$ decreases with the decrease of parameters $\gamma$ and $\alpha_r$.
Compared with the parameter $\alpha_r$, the parameter $\gamma$ has a more significant impact on the gap $g$.
When $\gamma = 1\times 10^{-6}$ and $\alpha_r = 0.1$, the gap $g$ is very small for different meshes, which indicates that the self-contact is well simulated.
Besides, the results for different meshes are very close, which indicates that the SFVEM has a good mesh convergence.

\begin{table}[htbp]
\centering
\caption{Gap $g$ between the surfaces of the upper and lower flange for different parameters and different meshes.}
\label{p10.s5.1.t1}
\begin{tabular}{lllll}
\hline
\multicolumn{1}{c}{$\alpha_r$} & \multicolumn{1}{c}{$\gamma$} & \multicolumn{1}{c}{$g$ (mesh 1)} & \multicolumn{1}{c}{$g$ (mesh 2)} & \multicolumn{1}{c}{$g$ (mesh 3)} \\ \hline
\multirow{3}{*}{10}            & $1\times 10^{-4}$            & 8.6882E-03                     & 8.7000E-03                     & 7.5000E-03                     \\
                               & $1\times 10^{-5}$            & 2.8000E-03                     & 1.9476E-03                     & 1.7330E-03                     \\
                               & $1\times 10^{-6}$            & 7.3795E-04                     & 8.2579E-04                     & 4.1061E-04                     \\ \hline
\multirow{3}{*}{1.0}           & $1\times 10^{-4}$            & 7.9000E-03                     & 7.2327E-03                     & 7.0869E-03                     \\
                               & $1\times 10^{-5}$            & 2.1822E-03                     & 2.1104E-03                     & 1.7138E-03                     \\
                               & $1\times 10^{-6}$            & 5.1084E-04                     & 5.3938E-04                     & 4.9040E-04                     \\ \hline
\multirow{3}{*}{0.1}           & $1\times 10^{-4}$            & 7.6000E-03                     & 7.3281E-03                     & 7.0544E-03                     \\
                               & $1\times 10^{-5}$            & 1.8730E-03                     & 1.8934E-03                     & 1.7435E-03                     \\
                               & $1\times 10^{-6}$            & 4.8469E-04                     & 4.8581E-04                     & 4.6363E-04                     \\ \hline
\end{tabular}
\end{table}

\begin{figure}[htbp]
    \centering
    \includegraphics[width=1.0\textwidth]{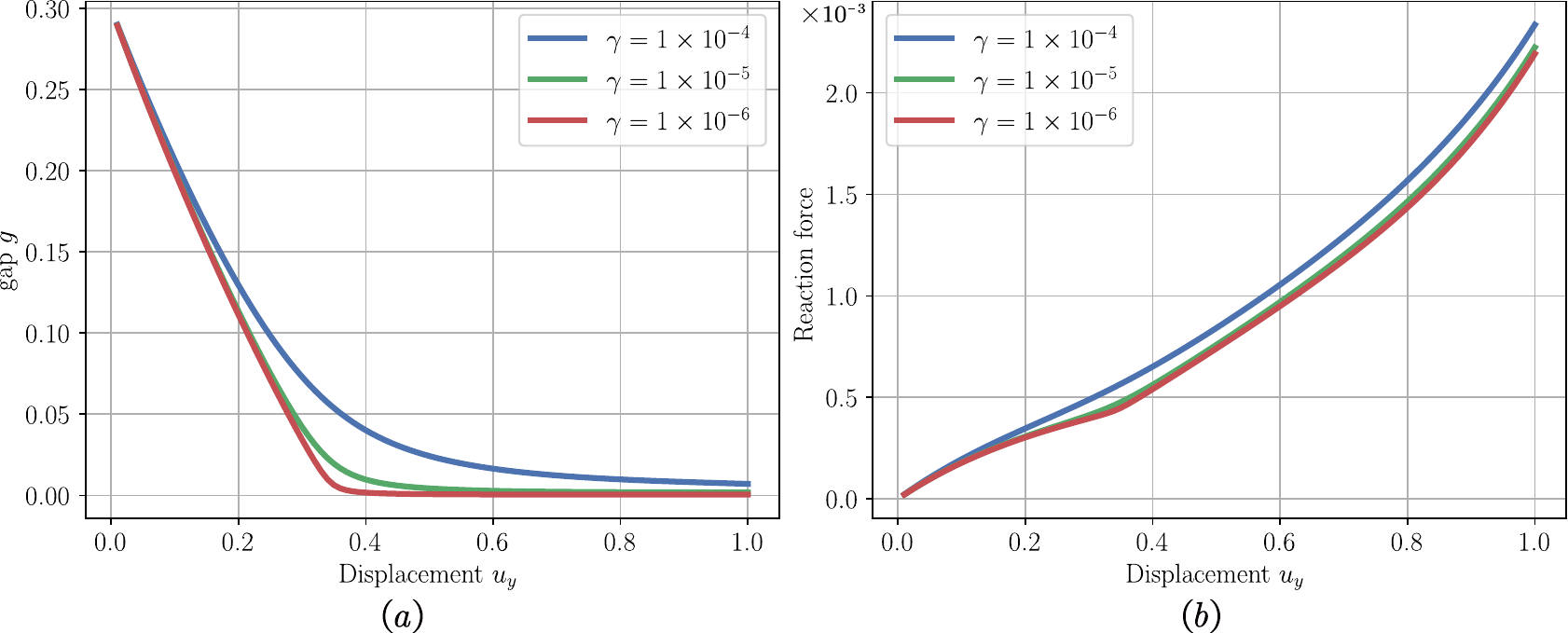}
    \caption{Gap $g$ between the surfaces of the upper and lower flange and reaction force during loading for different parameter $\gamma$, for mesh 3 and $\alpha_r = 0.1$.}
    \label{p10.s5.1.f2}
\end{figure}

In order to text the influence of the parameter $\gamma$ detailly,
we select $\alpha_r = 0.1$ for mesh 3 to investigate the influence of the parameter $\gamma$ on the gap function $g$ during loading.
Besides, the reaction force at the top of the box is also recorded.
As shown in Fig.\ref{p10.s5.1.f2} (a), the gap function $g$ is sensitive with the parameter $\gamma$. 
For large $\gamma$, the third medium can provide a large reaction force, 
resulting in poor gap closure. 
For small $\gamma$, the gap is nearly closed, so the third medium approach can approximate the contact behavior and is consistent with simulations using the classical contact discretization.
Similar result can also be found from the reaction force, see Fig.\ref{p10.s5.1.f2} (b).

Next, the influence of different regularization terms is investigated.
Here we consider three different regularization terms as given in Eq.\eqref{p10.s2.3.Psir1}, Eq.\eqref{p10.s2.3.rw1} and Eq.\eqref{p10.s2.3.rw2} (also see Table \ref{p10.s2.3.t1}).
According to the previous analysis, the parameters in the third medium can be selected as $\gamma = 1\times 10^{-6}$ and $\alpha_r = 10$.
For different regularization terms, 
the gap $g$ between the surfaces of the upper and lower flange is recorded in Table \ref{p10.s5.1.t2}.
It can be seen that the choice of regularization term has little effect on the results.
Similar results can also be obtained if the item $e^{-\beta |F|}$ is considered in the regularization term.

\begin{table}[htbp]
    \caption{Gap $g$ between the surfaces of the upper and lower flange for different regularization terms, for $\gamma = 1\times 10^{-6}$ and $\alpha_r = 10$.}
    \label{p10.s5.1.t2}
    \centering
    \begin{tabular}{lccc}
    \hline
    \multicolumn{1}{c}{$\Psi_r(\bm{u})$}                                                                                                                                                   & $g$ (mesh 1) & $g$ (mesh 2) & $g$ (mesh 3) \\ \hline
    $\frac{1}{2} \left[\nabla\bm{F}\vdots\nabla\bm{F} - \frac{1}{2} \text{Div}(\nabla\bm{u})\cdot\text{Div}(\nabla\bm{u})\right]$                                & 5.0282E-04   & 5.0988E-04   & 4.9311E-04   \\
    $\frac{1}{2} \left(\nabla\varphi\cdot\nabla\varphi+\nabla J\cdot\nabla J\right)$                                                                                   & 5.1084E-04   & 5.3938E-04   & 4.1061E-04   \\
    $\frac{1}{2} \left[\nabla\left[\frac{F_{12}-F_{21}}{F_{11}+F_{22}}\right]\cdot\nabla\left[\frac{F_{12}-F_{21}}{F_{11}+F_{22}}\right]+\nabla J\cdot\nabla J\right]$ & 5.1127E-04   & 5.4095E-04   & 4.0902E-04   \\ \hline
    \end{tabular}
\end{table}

\begin{figure}
    \centering
    \includegraphics[width=1.0\textwidth]{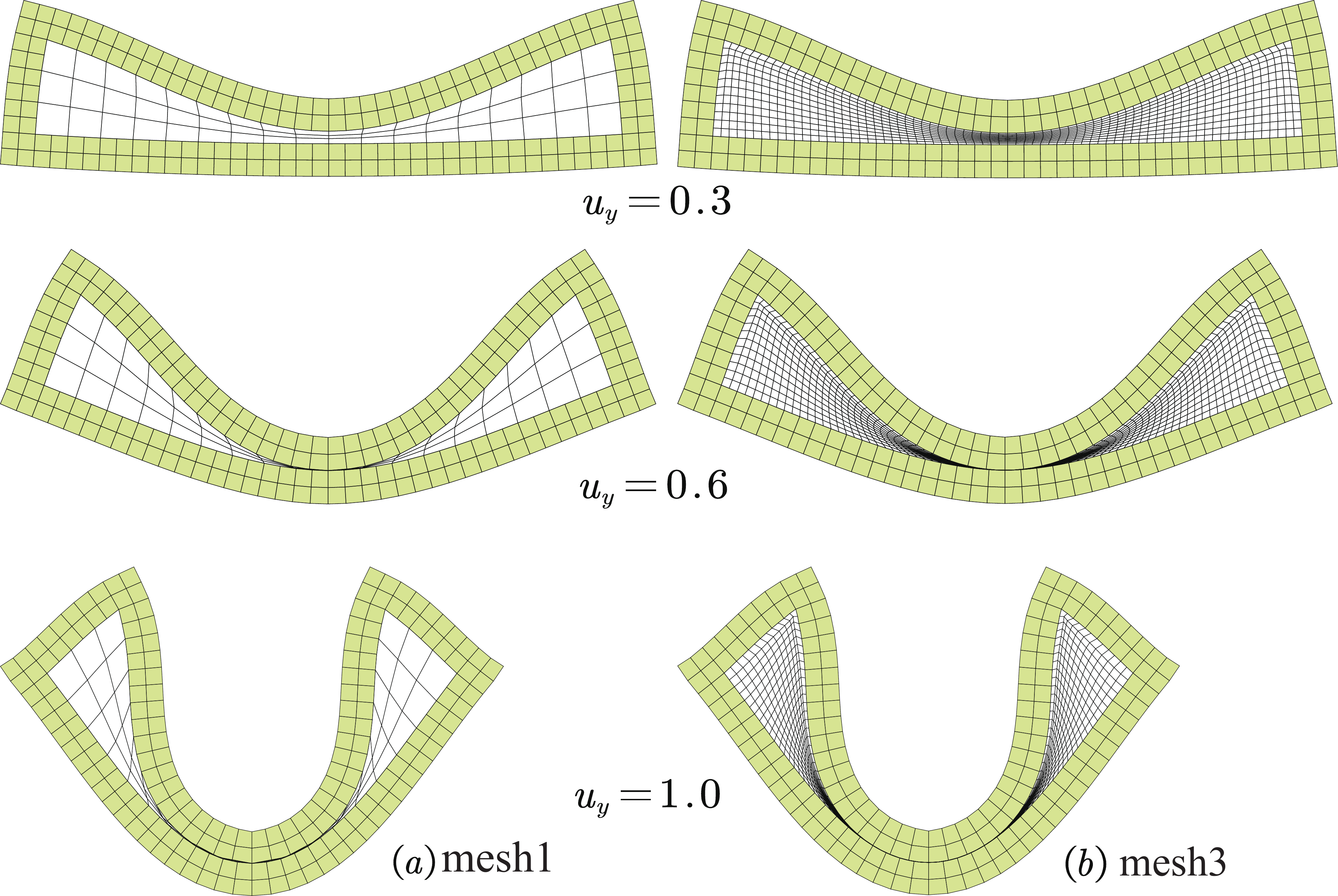}
    \caption{Deformed configuration for the box with self-contact under different displacement $u_y$, for mesh 1 and mesh 3, $\alpha_r = 0.1$ and $\gamma = 1\times 10^{-6}$.}
    \label{p10.s5.1.f3}
\end{figure}

For different meshes and different displacement $u_y$,
the deformed configuration of the box with self-contact is shown in Fig.\ref{p10.s5.1.f3}.
It can be seen that the SFVEM can well simulate the self-contact problem of the box with the third medium contact.
Besides, the results obtained by SFVEM attach well with the results obtained by FEM under the third medium contact framework, see \cite{TMC2} for detail.
Compared with FEM, SFVEM can use polygonal elements, which makes the mesh generation and local refinement easier.

\subsection{Self-contact of the C-box}
\label{p10.s5.2}
\begin{figure}[htbp]
    \centering
    \includegraphics[width=1.0\textwidth]{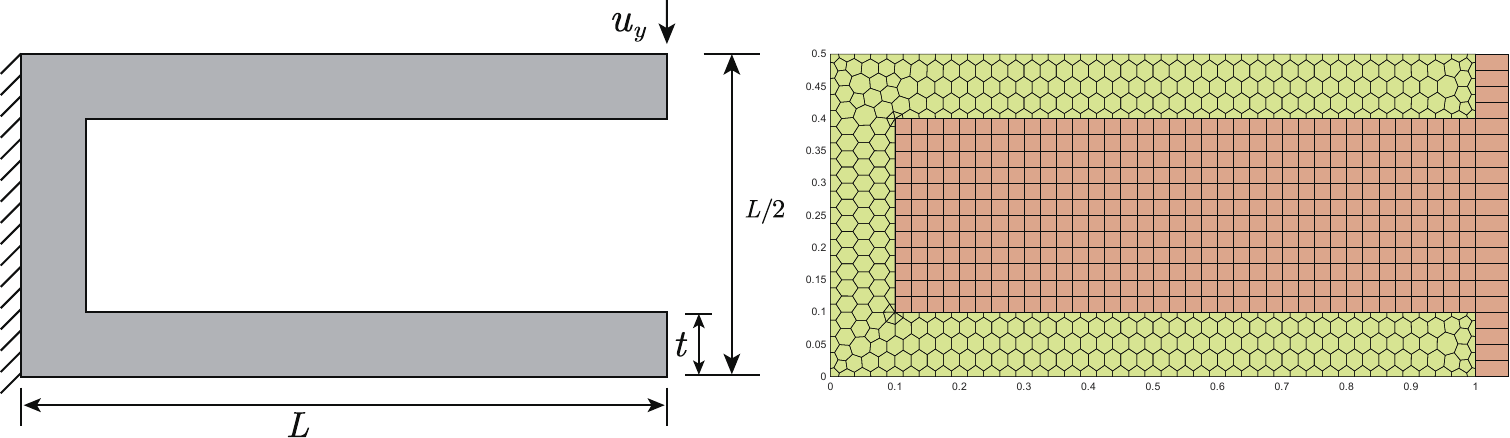}
    \caption{Self-contact of the C-box, geometry model, boundary condition and VEM mesh.}
    \label{p10.s5.2.f1}
\end{figure}

Another self-contact problem is the benchmark originally proposed in \cite{TMC3}.
The geometry of the C-box and boundary conditions are illustrated in Fig.\ref{p10.s5.2.f1}.
The geometry parameters are selected as $L = 1$ and $t = 0.1$. 
Under current assumption, the initial gap between the upper beam and the lower beam is $g_0 = 0.3$.
The left side of the C-box if fixed and the upper beam is loaded by a prescribed displacement $u_y$ at the right top point, as shown in Fig.\ref{p10.s5.2.f1}. 
The material parameters are chosen as $K=5/3$ and $\mu=5/14$ which yields very flexible beams.

\begin{figure}[htbp]
    \centering
    \includegraphics[width=1.0\textwidth]{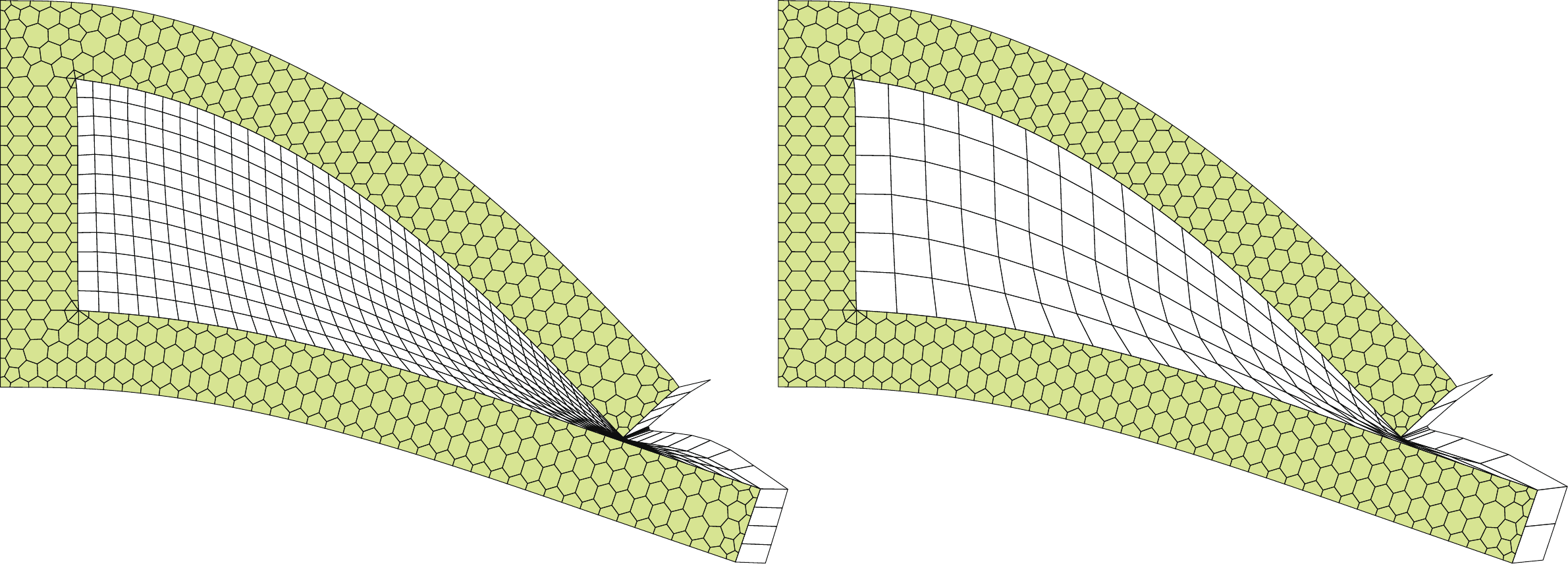}
    \caption{Deformed configuration for the C-box with self-contact under different meshes.
    $u_y = -0.5$, $\gamma = 1\times 10^{-5}$ and $\alpha_r = 1$}
    \label{p10.s5.2.f2}
\end{figure}

Under displacement load, the upper and lower parts of the structure will come into contact (point contact of the upper beam with the lower one).
As shown in Fig.\ref{p10.s5.2.f2}, the third medium is used in the entire cavity.
Besides, there is a column of third medium elements at the right side of the system to provide a sufficient number of elements in the vicinity of the contact point.
The polygonal mesh is used for the solid domain and quadrilateral mesh is used for the third medium.
In this example, 
the regularization term is selected as
$\Psi_r = \frac{1}{2} \left(\nabla\varphi\cdot\nabla\varphi+\nabla J\cdot\nabla J\right)$ as mentioned in Eq.\eqref{p10.s2.3.Psir1},
Then the potential energy for the third medium domain can be expressed as
\begin{equation}
    \label{p10.s5.2.WTMC}
    W^{TMC}(\bm{u}) = \int_{\Omega_m}\gamma \left[\frac{\mu}{2}\left(J^{-\frac{2}{3}}\text{tr}(\bm{C})-3\right)
    +\alpha_r \frac{1}{2} \left(\nabla\varphi\cdot\nabla\varphi+\nabla J\cdot\nabla J\right)\right]\ud\Omega.
\end{equation}

Firstly, we consider a vertical displacement $u_y = -0.5$.
The parameters in the third medium are selected as $\gamma = 1\times 10^{-5}$ and $\alpha_r = 1$.
Besides, for the third medium region, we consider two meshes with different densities.
Since non-matching polygonal elements can be used in SFVEM, no changes are required to the mesh of the solid domain.
The deformed states for two different meshes are shown in Fig.\ref{p10.s5.2.f2}.
The deformed mesh illustrates clearly the characteristic of a point contact for this loading assumption.

\begin{figure}[htbp]
    \centering
    \includegraphics[width=1.0\textwidth]{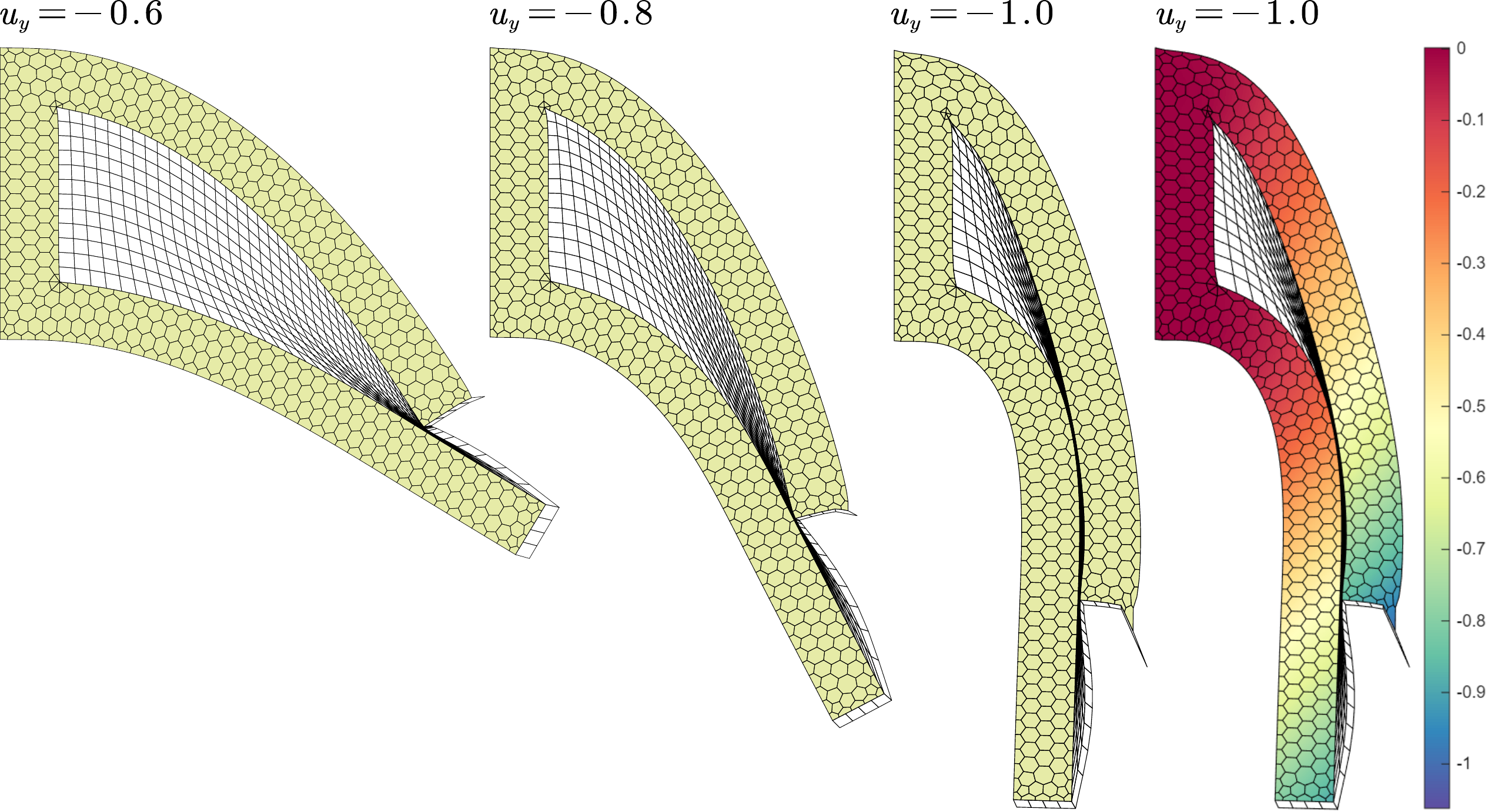}
    \caption{Deformed configuration for the C-box with self-contact under different displacement $u_y$ and contour plot of $u_y$.}
    \label{p10.s5.2.f3}
\end{figure}

Then we can consider a very large load with $u_y = -1$.
In order to ensure the convergence of the calculation, 
the parameters in the third medium need to be appropriately enlarged. 
Here, we choose $\gamma = 1\times 10^{-5}$ and $\alpha_r = 20$.
The deformed configuration for the C-box under different displacement $u_y$ and contour plot of displacement $u_y$ at the last state are shown in Fig.\ref{p10.s5.2.f3}.
Easy to find that the mesh in the extremely deformed third medium is controlled well by the regularization term.

\subsection{Punch problem}
\label{p10.s5.3}
\begin{figure}[htbp]
    \centering
    \includegraphics[width=0.6\textwidth]{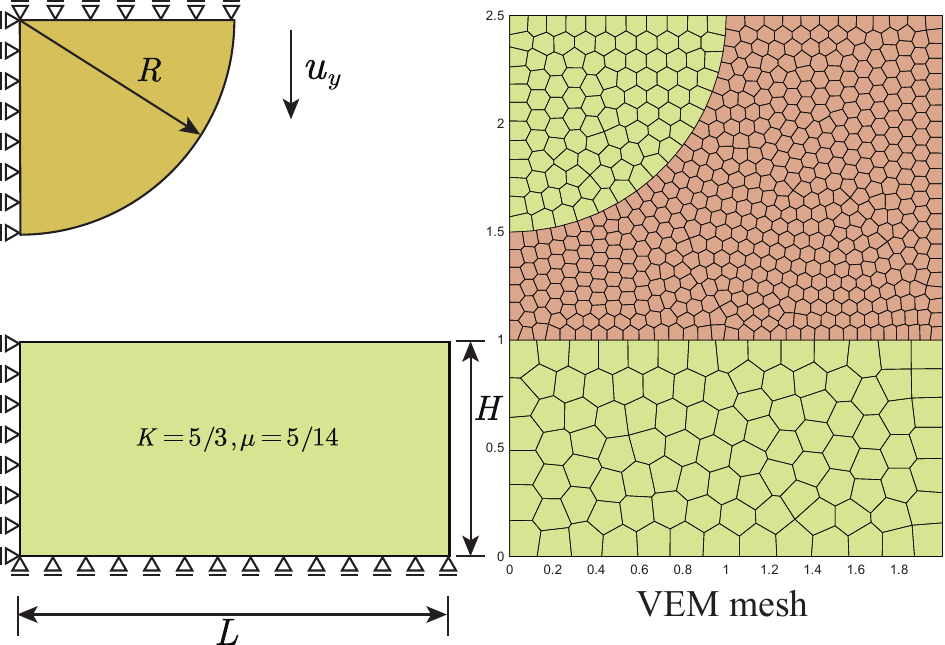}
    \caption{Punch problem, geometry model and boundary conditions.}
    \label{p10.s5.3.f1}
\end{figure}

In this example, a classic punch problem is considered.
The contact problem consists of two parts: the lower rectangular pressure domain and the upper pressure domain. 
Due to symmetry, only half is considered in the calculation.
As given in Fig.\ref{p10.s5.3.f1}, the size of the rectangle is $L\times H = 2\times 1$ and the radius of the circle is $R = 1$.
The displacement boundary conditions are also shown in Fig.\ref{p10.s5.3.f1}.
A vertical displacement load $u_y$ is applied at the top of the upper domain.
The material parameter for the rectangular is selected as $K=5/3$ and $\mu=5/14$ which yields very flexible bodies.
Beisdes, the different material parameters for the circular punch are selected for comparison.

The third medium is used in the entire cavity between the two bodies.
As shown in Fig.\ref{p10.s5.3.f1}, the polygonal meshes are used for the solid domain and the third medium.
The regularization term is selected as
$\Psi_r = \frac{1}{2} \left(\nabla\varphi\cdot\nabla\varphi+\nabla J\cdot\nabla J\right)$, which is the same as Eq.\eqref{p10.s2.3.Psir1}.
Then the potential energy for the third medium domain can be expressed as Eq.\eqref{p10.s5.2.WTMC}.
For this problem, since the third medium is compressed and deformed very greatly, 
extremely small parameters will lead to non-convergence of the calculation. 
Here the parameters are selected as $\gamma = 1\times 10^{-4}$ and $\alpha = 1$.

\begin{figure}[htbp]
    \centering
    \includegraphics[width=1.0\textwidth]{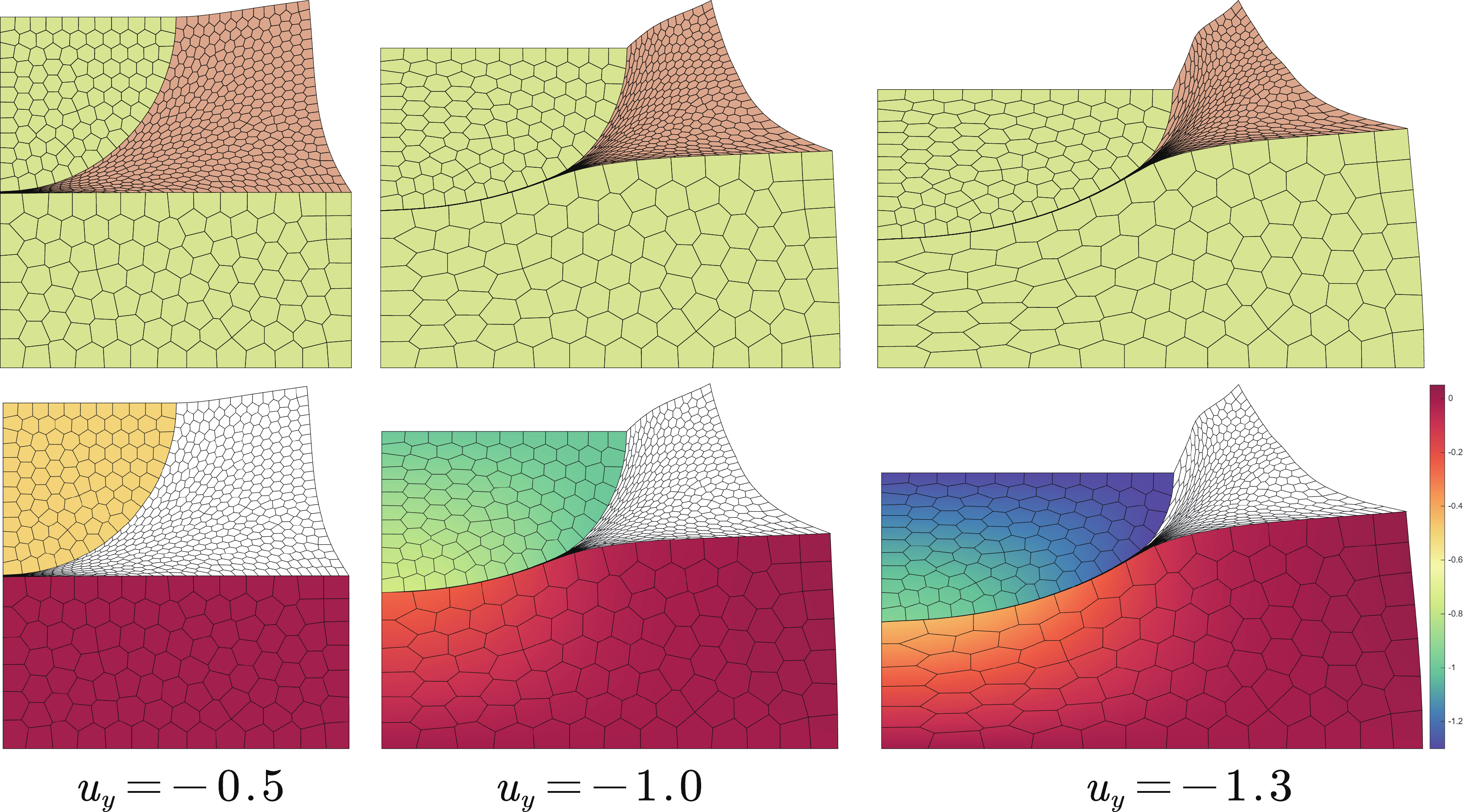}
    \caption{Deformed configuration and contour plots of $u_y$ under different displacement loads, material parameters for the circular punch are selected as $K=5/3$ and $\mu=5/14$.}
    \label{p10.s5.3.f2}
\end{figure}

We consider a vertical displacement $u_y = -1.3$ and $130$ loading steps within the Newton-Raphson solution scheme.
The material parameters for the circular punch are selected as $K=5/3$ and $\mu=5/14$ (same as the rectangular domain).
Under this condition, the deformed configuration and contour plot of displacement $u_y$ for different times are shown in Fig.\ref{p10.s5.3.f2}.
The mesh in the third medium domain is severely compressed. 
The existence of regularization makes the convergence of the calculation well controlled.

\begin{figure}[htbp]
    \centering
    \includegraphics[width=1.0\textwidth]{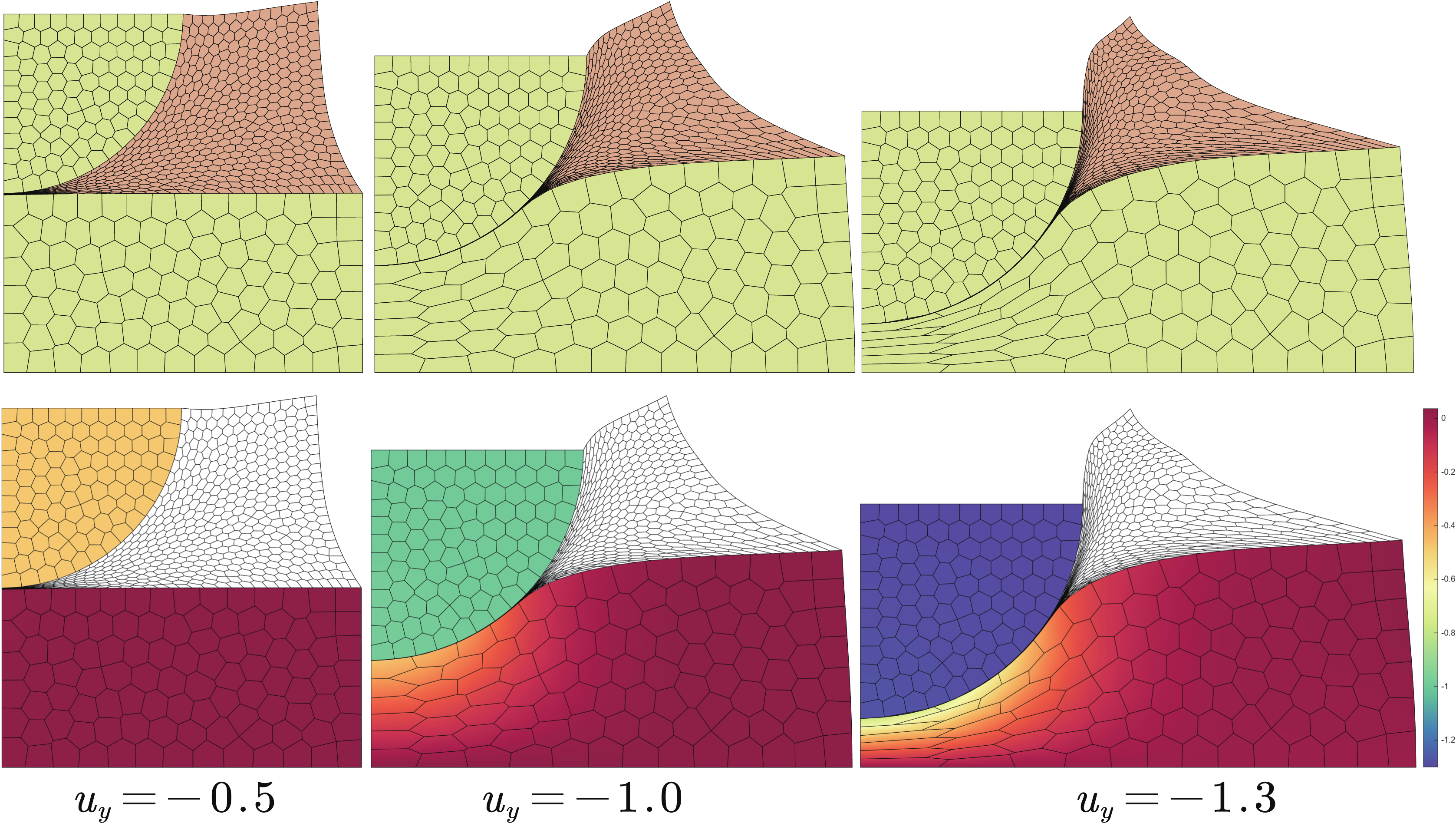}
    \caption{Deformed configuration and contour plots of $u_y$ under different displacement loads, material parameters for the circular punch are selected as $K=500/3$ and $\mu=500/14$.}
    \label{p10.s5.3.f3}
\end{figure}

Then we consider different material parameters as $K=500/3$ and $\mu=500/14$ for the circular punch, which yields a very rigid punch.
The deformed configuration and contour plot of displacement $u_y$ for different times are also shown in Fig.\ref{p10.s5.3.f3}.
It can be seen that the lower rectangle has undergone a larger deformation due to its softness.
In this case, convergence of the problem is more difficult than for the first material assumption ($K=5/3$ and $\mu=5/14$ for the circular punch), 
so an automatic load step adjustment is required.

\subsection{Complex contact of multiple objects}
\label{p10.s5.4}
\begin{figure}[htbp]
    \centering
    \includegraphics[width=1.0\textwidth]{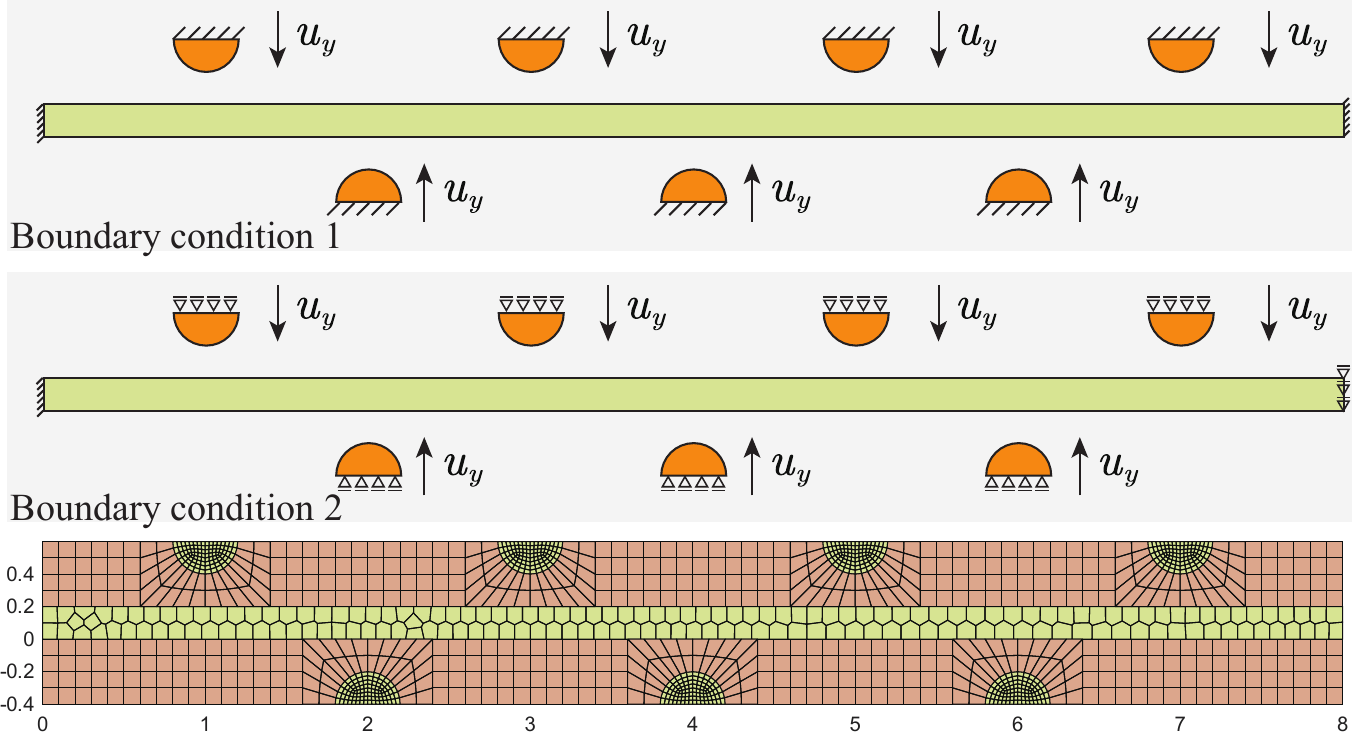}
    \caption{Geometry model and different boundary conditions.}
    \label{p10.s5.4.f1}
\end{figure}

In this example, we consider a complex contact problem of multiple objects.
As shown in Fig.\ref{p10.s5.4.f1}, 
the geometry consists of different parts: 
a rectangle is fixed and 7 semicircles are placed around the rectangle.
The size of the rectangle is $L\times H = 8\times 0.2$ and the radius of the semicircle is $R = 0.2$.
The third medium is used in the entire cavity between the different parts.
The polygonal meshes are used for the solid domain, as shown in Fig.~\ref{p10.s5.4.f1}.
The material parameters for the rectangular part are selected as $K=50$ and $\mu=10$.
At the same time, the material parameters for the semicircular parts are selected as $K=50000$ and $\mu=10000$ 
which yields very hard bodies.

\begin{figure}[htbp]
    \centering
    \includegraphics[width=1.0\textwidth]{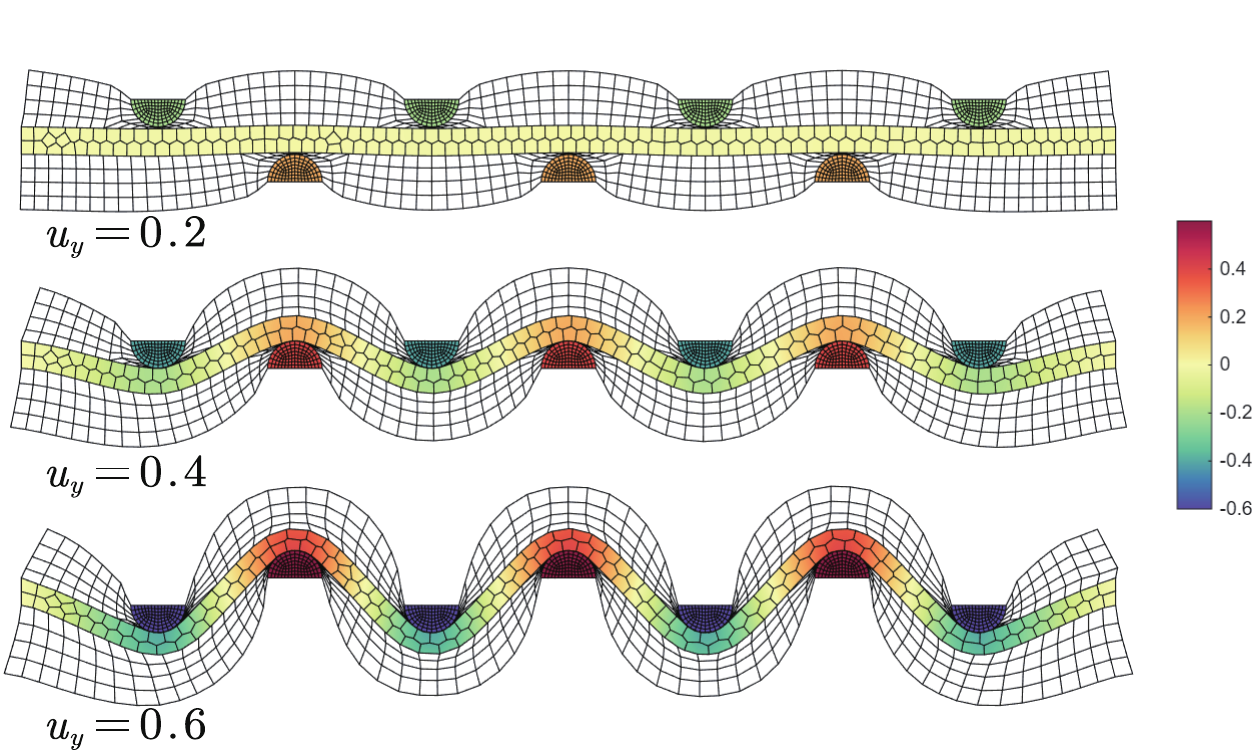}
    \caption{Deformed configuration and contour plots of $u_y$ under different displacement loads, for the first boundary condition.}
    \label{p10.s5.4.f2}
\end{figure}

The regularization term for the third medium is selected as
$\Psi_r = \frac{1}{2} \left(\nabla\varphi\cdot\nabla\varphi+\nabla J\cdot\nabla J\right)$, which is the same as Eq.\eqref{p10.s2.3.Psir1}.
Then the potential energy for the third medium domain can be expressed as Eq.\eqref{p10.s5.2.WTMC}.
The parameters are selected as $\gamma = 1\times 10^{-4}$ and $\alpha = 10$ for different boundary conditions.
In this example,
two different boundary conditions are considered.
For the first boundary condition,
the rectangle is fixed at the both ends, 
while the seven middle semicircles are subjected to a prescribed vertical displacement $u_y$ with their motion constrained in the $x$-direction (Fig.~\ref{p10.s5.4.f1}, top).
Under this condition, the deformed configuration and contour plot of displacement $u_y$ for different times are shown in Fig.\ref{p10.s5.4.f2}.

\begin{figure}[htbp]
    \centering
    \includegraphics[width=1.0\textwidth]{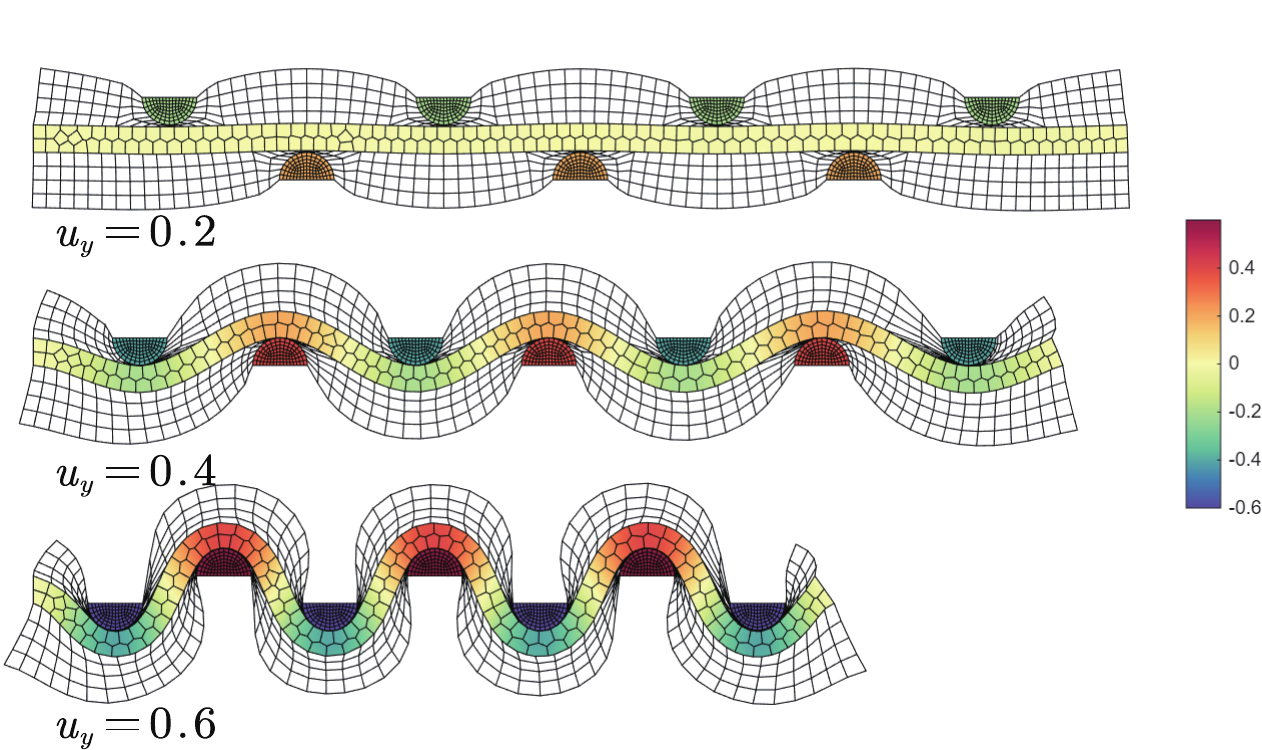}
    \caption{Deformed configuration and contour plots of $u_y$ under different displacement loads, for the second boundary condition.}
    \label{p10.s5.4.f3}
\end{figure}

Then, the rectangle is fixed at the left end and simply supported at the right end,
while the seven middle semicircles are subjected to a prescribed vertical displacement $u_y$ (Fig.~\ref{p10.s5.4.f1}, middle).
In this case, the seven semicircular components will have displacement in the $x$-direction as the displacement load is applied, 
which will cause the sliding of the third medium mesh.
Under this condition, the deformed configuration and contour plot of displacement $u_y$ for different times are shown in Fig.\ref{p10.s5.4.f3}.
Compared with the first boundary condition, the seven middle semicircles have the displacement in the $x$-direction,
which leads to a more complex contact state.
Easy to find that the SFVEM can well simulate the complex contact problem of multiple objects with the third medium contact.
Obviously, the third medium contact avoids the complex contact surface search and treatment of inequality constraints.

\section{Conclusion}
\label{p10.s6}
In this work, we have developed a stabilization-free virtual element method for the third medium contact in 2D problems.
The SFVEM is used to discretize both the solid bodies and the third medium.
Since SFVEM does not require any additional stabilization terms, it has more advantages in dealing with nonlinear problems.
Compared with the traditional contact algorithms,
the third medium contact avoids the complex contact surface search and treatment of inequality constraints.
The discretization form of the third medium contact is given in detail.
Furthermore, different regularization terms for the third medium are also discussed and compared in the numerical examples.
Several numerical examples are presented to verify the effectiveness of the SFVEM for the third medium contact.
The numerical results show that the SFVEM can well simulate the contact problems under different situations.
The first-order VEMs can be used by combining with the techniques in \cite{TMC2,TMCWriggers2}.
Besides, we will extend the current work to 3D problems and frictional contact problems in the next work.

\appendix
\section{Constitutive and stress tensor for the hyperelastic body}
\label{p10.appendix1}
\lstset{
    basicstyle=\footnotesize\ttfamily\setstretch{1}, 
    columns=fixed,       
    numbers=left,                                        
    numberstyle=\tiny\color{gray},                       
    frame=none,                                          
    backgroundcolor=\color[RGB]{245,245,244},            
    keywordstyle=\color[RGB]{40,40,255},                 
    numberstyle=\footnotesize\color{darkgray},           
    commentstyle=\it\color[RGB]{0,96,96},                
    stringstyle=\rmfamily\slshape\color[RGB]{128,0,0},   
    showstringspaces=false,                              
    language=matlab                                        
}

\begin{lstlisting}
    clear;
    syms K mu
    C = sym('C', [3,3],'real'); % Right Cauchy-Green tensor
    I1 = trace(C);   
    J = sqrt(det(C));
    W = K/2*(log(J))^2+mu/2*(J^(-2/3)*I1-3);
    % Derivation of second Piola-Kirchhoff stress
    S = sym(zeros(3,3));
    for i = 1:3
        for j = 1:3
            S(i,j) = 2*diff(W,C(i,j));
        end
    end
    % Derivation of constitutive tensor
    D = sym(zeros(3,3,3,3));
    for i = 1:3
        for j = 1:3
            for k = 1:3
                for l = 1:3
                    D(i,j,k,l) = 2*diff(S(i,j),C(k,l));
                end
            end
        end
    end
\end{lstlisting}

\section{Constitutive and stress tensor for the third medium}
\label{p10.appendix2}
\begin{lstlisting}
clear;
syms beta;
% deformation gradient and its gradient
F = sym('F', [2,2],'real');
dF = sym('dF', [2,2,2],'real');
% regularization term,
Lu = [dF(1,1,1)+dF(1,2,2);dF(2,1,1)+dF(2,2,2)]; % Div(\nabla u)
Psi = 1/2*(dF(:)'*dF(:)-1/2*(Lu'*Lu))*exp(-beta*det(F));
P = sym(zeros(2,2));
T = sym(zeros(2,2,2));
for i = 1:2
    for j = 1:2
        P(i,j) = diff(Psi,F(i,j));
    end
end
for i = 1:2
    for j = 1:2
        for k = 1:2
            T(i,j,k) = diff(Psi,dF(i,j,k));
        end
    end
end
A = sym(zeros(2,2,2,2,2));
for i = 1:2
    for j = 1:2
        for k = 1:2
            for m = 1:2
                for n = 1:2
                    A(i,j,k,m,n) = diff(T(i,j,k),F(m,n));
                end
            end
        end
    end
end
B = sym(zeros(2,2,2,2,2,2));
for i = 1:2
    for j = 1:2
        for k = 1:2
            for m = 1:2
                for n = 1:2
                    for p = 1:2
                        B(i,j,k,m,n,p) = diff(T(i,j,k),dF(m,n,p));
                    end
                end
            end
        end
    end
end

D = sym(zeros(2,2,2,2));
for i = 1:2
    for j = 1:2
        for m = 1:2
            for n = 1:2
                D(i,j,m,n) = diff(P(i,j),F(m,n));
            end
        end
    end
end
\end{lstlisting}

\section{Matrix forms of the stress and constitutive tensors for the third medium}
\label{p10.appendix3}
The matrix forms of the stress and constitutive tensors in Eq.\eqref{p10.s4.2.ddw} and Eq.\eqref{p10.s4.2.Ke} for the third medium are given as follows.
\begin{equation}
	\hat{\mathbb{D}} = \begin{bmatrix}
		\mathbb{D}_{1111} & \mathbb{D}_{1121} & \mathbb{D}_{1112} & \mathbb{D}_{1122}\\
		\mathbb{D}_{2111} & \mathbb{D}_{2121} & \mathbb{D}_{2112} & \mathbb{D}_{2122}\\
		\mathbb{D}_{1211} & \mathbb{D}_{1221} & \mathbb{D}_{1212} & \mathbb{D}_{1222}\\
		\mathbb{D}_{2211} & \mathbb{D}_{2221} & \mathbb{D}_{2212} & \mathbb{D}_{2222}\\
	\end{bmatrix},
\end{equation}
\begin{equation}
	\hat{\mathbb{A}} = \begin{bmatrix}
		\mathbb{A}_{11111} & \mathbb{A}_{11121} & \mathbb{A}_{11112} & \mathbb{A}_{11122}\\
		\mathbb{A}_{21111} & \mathbb{A}_{21121} & \mathbb{A}_{21112} & \mathbb{A}_{21122}\\
		\mathbb{A}_{12111} & \mathbb{A}_{12121} & \mathbb{A}_{12112} & \mathbb{A}_{12122}\\
		\mathbb{A}_{22111} & \mathbb{A}_{22121} & \mathbb{A}_{22112} & \mathbb{A}_{22122}\\
		\mathbb{A}_{11211} & \mathbb{A}_{11221} & \mathbb{A}_{11212} & \mathbb{A}_{11222}\\
		\mathbb{A}_{21211} & \mathbb{A}_{21221} & \mathbb{A}_{21212} & \mathbb{A}_{21222}\\
		\mathbb{A}_{12211} & \mathbb{A}_{12221} & \mathbb{A}_{12212} & \mathbb{A}_{12222}\\
		\mathbb{A}_{22211} & \mathbb{A}_{22221} & \mathbb{A}_{22212} & \mathbb{A}_{22222}\\
	\end{bmatrix},
\end{equation}
and 
\begin{equation}
	\hat{\mathbb{B}} = \begin{bmatrix}
		\mathbb{B}_{111111} & \mathbb{B}_{111211} & \mathbb{B}_{111121} & \mathbb{B}_{111221}& \mathbb{B}_{111112} & \mathbb{B}_{111212} & \mathbb{B}_{111122} & \mathbb{B}_{111222}\\
		\mathbb{B}_{211111} & \mathbb{B}_{211211} & \mathbb{B}_{211121} & \mathbb{B}_{211221}& \mathbb{B}_{211112} & \mathbb{B}_{211212} & \mathbb{B}_{211122} & \mathbb{B}_{211222}\\
		\mathbb{B}_{121111} & \mathbb{B}_{121211} & \mathbb{B}_{121121} & \mathbb{B}_{121221}& \mathbb{B}_{121112} & \mathbb{B}_{121212} & \mathbb{B}_{121122} & \mathbb{B}_{121222}\\
		\mathbb{B}_{221111} & \mathbb{B}_{221211} & \mathbb{B}_{221121} & \mathbb{B}_{221221}& \mathbb{B}_{221112} & \mathbb{B}_{221212} & \mathbb{B}_{221122} & \mathbb{B}_{221222}\\
		\mathbb{B}_{112111} & \mathbb{B}_{112211} & \mathbb{B}_{112121} & \mathbb{B}_{112221}& \mathbb{B}_{112112} & \mathbb{B}_{112212} & \mathbb{B}_{112122} & \mathbb{B}_{112222}\\
		\mathbb{B}_{212111} & \mathbb{B}_{212211} & \mathbb{B}_{212121} & \mathbb{B}_{212221}& \mathbb{B}_{212112} & \mathbb{B}_{212212} & \mathbb{B}_{212122} & \mathbb{B}_{212222}\\
		\mathbb{B}_{122111} & \mathbb{B}_{122211} & \mathbb{B}_{122121} & \mathbb{B}_{122221}& \mathbb{B}_{122112} & \mathbb{B}_{122212} & \mathbb{B}_{122122} & \mathbb{B}_{122222}\\
		\mathbb{B}_{222111} & \mathbb{B}_{222211} & \mathbb{B}_{222121} & \mathbb{B}_{222221}& \mathbb{B}_{222112} & \mathbb{B}_{222212} & \mathbb{B}_{222122} & \mathbb{B}_{222222}\\
	\end{bmatrix}.
\end{equation}


\bibliography{contact}

\end{document}